\documentclass[reqno,11pt, oneside]{amsart}
\usepackage{enumerate}
\usepackage[nobysame]{amsrefs}
\usepackage{enumitem}
\usepackage{amssymb}
\usepackage{amsmath}
\usepackage{amsthm}
\usepackage{microtype}
\usepackage{graphicx}
\usepackage{array}
\usepackage{fullpage}
	
\theoremstyle{plain}
\newtheorem{theorem}{Theorem}[section]
\newtheorem{corollary}[theorem]{Corollary}
\newtheorem{lemma}[theorem]{Lemma}
\newtheorem{proposition}[theorem]{Proposition}

\theoremstyle{definition}
\newtheorem{remark}[theorem]{Remark}
\newtheorem{example}[theorem]{Example}

\newtheorem{definition}[theorem]{Definition}

\numberwithin{equation}{section}

\newcommand{\Z}{\in \mathbb Z}

\newcommand{\f}{f=\sum\limits_{n\in\mathbb Z}a_nz^n}
\newcommand{\g}{g=\sum\limits_{n\in\mathbb Z}b_nz^n}

\title{The inverse and the composition in the set of formal Laurent series}
\author{Dawid Bugajewski$^1$ \\F\lowercase{aculty of }P\lowercase{hysics, }U\lowercase{niversity of }W\lowercase{arsaw, }W\lowercase{arsaw, }P\lowercase{oland}}\thanks{$^1$E-mail address: d.bugajewski@student.uw.edu.pl; ORCID number: 0009-0002-2327-6037}


\begin{document}

\begin{abstract}The aim of this article is to investigate the issues of multiplicative inverses and composition in the set of formal Laurent series. We show the lack of general uniqueness of inverses of formal Laurent series; necessary and sufficient conditions for the existence of inverses of formal Laurent series satisfying some weak assumptions are also provided. Moreover, we define a general composition of formal Laurent series and investigate the Right Distributive Law and the Chain Rule in this context. 
\end{abstract}

\maketitle

\section{Introduction}\label{sec1}
A formal power series over a ring $S$ is defined as a mapping $f:(\mathbb Z_+\cup\{0\})\rightarrow S$; we write $f=\sum\limits_{n=0}^{\infty}a_nz^n$, where $a_n:=f(n)$. The set of all formal power series over $S$ is usually denoted as either $S[[z]]$ (see e.g. \cite{radicals}) or $\mathbb X(S)$ (see e.g. \cite{zrodlonew}). It is well-known that some of the properties of $\mathbb X(S)$ depend strongly on what ring $S$ is considered -- for instance, it follows simply from the definition of multiplication of formal power series $\left(\sum\limits_{n=0}^{\infty}a_nz^n\right)\left(\sum\limits_{n=0}^{\infty}b_nz^n\right):=\sum\limits_{n=0}^{\infty}\left(\sum\limits_{k=0}^na_kb_{n-k}\right)z^n$ that if $S$ is an integral domain, then so is $\mathbb X(S)$; throughout this article, we will focus on the case $S=\mathbb C$. 

For a more systematic description of the structure of $\mathbb X(\mathbb C)$ (or, more generally, $\mathbb X(S)$), the reader is referred to e.g. \cite{zrodlo12,zrodlo3,zrodlonew}. Here let us emphasize that over decades, this topic has found various applications in, among others, algebra (see e.g. \cite{zrodlo15,zrodlo16,jcp}), differential equations (see e.g. \cite{zrodlo11,zrodlo17}), number theory (see e.g. \cite{zrodlo12}) and combinatorics (see e.g. \cite{riordanbook}). Also, an important step discussed in more details in the Preliminaries was the development of the theory of composition of nonunit ($f(0)=0$) formal power series (see e.g. \cite{jennings}) and, at the beginning of the 21st century, its extension to all formal power series, a center result of which is the so called General Composition Theorem -- the necessary and sufficient condition for the existence of the general composition of any two formal power series proved by Gan and Knox in 2002 \cite{zrodlo8}. Let us add that some later, simpler proofs of that result can be found in e.g. \cite{zrodlo14} and \cite{zrodlonew}.

Formal Laurent series seem to be a natural generalization of formal power series, allowing for terms of the type $z^{-n}$, $n>0$. Although originally only a finite amount of such terms was allowed (a series satisfying this condition are now sometimes called formal semi-Laurent series, see e.g. \cite{zrodlonew}), the general case is also a subject of interest \cite{zrodlo5}. In this paper, we will focus on the latter (that is allow an arbitrary amount of nonzero terms with negative exponents), focusing on the space of formal Laurent series over $\mathbb C$ like \cite{zrodlo5}. One should point out that their properties are quite different from the analogous properties of formal power series -- for example, in general, the product of two given formal Laurent series may not exist \cite{zrodlo5} and the inverse in the set of formal Laurent series does not have to be unique, unlike to formal power series (see Examples in Section 3 of this article). Nevertheless, it should be emphasized that results concerning this general setting, despite some peculiar properties, apply to various important classes of formal Laurent series like for instance formal power series, formal semi-Laurent series, reversed formal semi-Laurent series (see e.g. \cite{zrodlonew}), the discrete group algebra $\ell^1(\mathbb Z)$, also denoted as $\mathbb L_1$, analogous sets $\mathbb L_p$ for $p>1$ or $p=+\infty$ \cite{zrodlo5} and others. Let us also add that the aforementioned paper \cite{zrodlo5} also introduces the concept of composition of formal Laurent series with formal power series. 

In this paper, we focus on two issues, namely (1) the multiplicative inverses in the set of all formal Laurent series over $\mathbb C$ and (2) extending the concept of composition of two formal series to the case of two formal Laurent series (or, in particular, a formal power series on the "outside" and a Laurent series on the "inside"). More precisely, in Section 3, we explore basic properties of multiplicative inverses of formal Laurent series -- we show that a formal Laurent series may have none, one or uncountably many inverses, and link this issue to the lack of associativity of multiplication of formal Laurent series and the existence of zero divisors. In Section 4, we address this topic more specifically, providing necessary and sufficient conditions for formal Laurent series satisfying some weak assumptions to be invertible and presenting some methods of calculating their inverses, if they exist. For that purpose we use the theory of infinite systems of linear algebraic equations, investigated for instance in \cite{determinant,periodic,zrodlo1,zrodlo13,fedorov2022,zrodlo2}. We present a few applications of the obtained methods -- for example we find all formal Laurent series inverses of unit formal power series (it is well-known that a formal power series $f$ satisfying $f(0)\neq 0$ possesses exactly one inverse in the set of formal power series, see e.g. \cite{zrodlonew}; here we find all its possible inverses, possibly including terms with negative exponents). Finally, in Section 5, we define a general composition of formal Laurent series and use it to examine the Right Distributive Law and the Chain Rule, which hold for  formal power series \cite{zrodlo9}.

\section{Preliminaries}\label{sec2}

In this section, we collect some basic definitions and results which will be needed in the sequel. Let us begin with the following simple
\begin{definition}\label{sums}
For a sequence $(a_n)_{n\in\mathbb Z}$ of complex numbers, we define the infinite sum $\sum\limits_{n\in\mathbb Z}a_n$ (or, equivalently, $\sum\limits_{n=-\infty}^{\infty}a_n$) as $\sum\limits_{n=-\infty}^{\infty}a_n:=\sum\limits_{n=1}^{\infty}a_{-n}+\sum\limits_{n=0}^{\infty}a_n$, provided both series on the right side are convergent. In other cases we say that the sum on the left side is divergent.
\end{definition}
\begin{remark}\label{rearr}
In the sequel, we will sometimes use the so called Rearrangement Theorem -- that is, for a function $a:\mathbb Z^k\rightarrow \mathbb C$ ($k\in\mathbb Z_+$, $k>1$) we will use the fact that if for a permutation $\sigma$ of $\{1,...,k\}$, the sum $\sum\limits_{n_{\sigma(1)}\in\mathbb Z}...\sum\limits_{n_{\sigma(k)}\in\mathbb Z}|a_{n_1,...,n_k}|$ converges, then for every permutation $\rho$ of $\{1,...,k\}$, $\sum\limits_{n_{\rho(1)}\in\mathbb Z}...\sum\limits_{n_{\rho(k)}\in\mathbb Z}a_{n_1,...,n_k}$ converges and the value of this sum does not depend on the choice of $\rho$. This is a simple conclusion from the Fubini and Tonelli theorems for $\mathbb Z^k$ with the product of $k$ counting measures; the additivity of the Lebesgue integral for finite measurable partitions ensures compatibility of this approach with Def. \ref{sums}. Let us, however, point out that there is also another way to prove this fact without using measure theory -- see e.g. \cite{zrodlo6}, Thm. 7.50, which proof can be easily extended to series over $\mathbb Z$ instead of $\mathbb N$ in the sense of Definition \ref{sums}.
\end{remark}
Let us now move on to the definition and basic properties of formal Laurent series.
\begin{definition}(see e.g. \cite{zrodlo5,zrodlonew,zrodlo3})
A formal Laurent series over $\mathbb C$ is defined as a mapping $g:\mathbb Z\rightarrow\mathbb C$. We denote $g=\sum\limits_{n\in\mathbb Z}b_nz^n$, where $b_n:=g(n)\in\mathbb C$ for every $n\in\mathbb Z$. The series $g^+:=\sum\limits_{n=0}^{\infty}b_nz^n$ and $g^-:=\sum\limits_{n=1}^{\infty}b_{-n}z^{-n}$ are called the regular and the principal part of $g$, respectively. We also denote $\breve{g}=\sum\limits_{n\in\mathbb Z}b_{-n}z^n$. If $b_n=0$ for all $n\in\mathbb Z$, we call $g$ the zero formal Laurent series and write simply $g=0$; analogously, if $b_n=\delta_{n,0}$ ($b_0=1$ and $b_n=0$ for $n\neq 0$), we call $g$ the unit formal Laurent series and write $g=1$. We denote by $\mathbb L (\mathbb C)$ (or simply by $\mathbb L$) the set of formal Laurent series over $\mathbb C$. \\
Now, let $\g,\f\in\mathbb L$, $c\in\mathbb C$. We define
\begin{enumerate}
\item the sum of $g$ and $f$: $g+f:=\sum\limits_{n\in\mathbb Z}(b_n+a_n)z^n$,
\item the product of $c$ and $g$ (scalar multiplication): $cg:=\sum\limits_{n\in\mathbb Z}(cb_n)z^n$,
\item the formal derivative of $g$: $g':=\sum\limits_{n\in\mathbb Z}(n+1)b_{n+1}z^n$,
\item the product of $g$ and $f$: $gf:=\sum\limits_{n\in\mathbb Z}\left(\sum\limits_{m\in\mathbb Z}b_ma_{n-m}\right)z^n$, provided all the sums over $m$ are convergent (otherwise we say $gf$ does not exist); we also denote $\mathbb L(g):=\{f\in\mathbb L:gf\mbox{ exists}\}$,
\item the natural powers of $g$: $g^0:=1$, $g^1:=1$, $g^k:=g^{k-1}g$ for $k>1$ (provided $g^{k-1}$ exists and $g\in\mathbb L(g^{k-1})$); we denote $g^k=\sum\limits_{n\in\mathbb Z}b_n^{(k)}z^n$ if it exists.
\end{enumerate}
\end{definition}
\begin{remark}\label{dlapopr}
It is obvious that if $g^k$ exists ($\g\in\mathbb L$, $k\in\mathbb N$, $k>1$), then 
\begin{eqnarray}\label{dopowers}
b_n^{(k)}=\sum\limits_{m_{k-1}\in\mathbb Z}\sum\limits_{m_{k-2}\in\mathbb Z}...\sum\limits_{m_{1}\in\mathbb Z}b_{n-m_{k-1}}b_{m_{k-1}-m_{k-2}}...b_{m_2-m_1}b_{m_1}.
\end{eqnarray}
Conversely, if the above series is convergent for all $n\in\mathbb Z$, then all powers of $g$ up to $g^k$ exist. Indeed, assume $g\neq 0$ and fix $m_{k-1}\in\mathbb Z$; then $b_{n-m_{k-1}}\neq 0$ for some $n\in\mathbb Z$. Therefore $\sum\limits_{m_{k-2}\in\mathbb Z}...\sum\limits_{m_{1}\in\mathbb Z}b_{m_{k-1}-m_{k-2}}...b_{m_2-m_1}b_{m_1}$ converges for all $m_{k-1}\in\mathbb Z$. Repeating this reasoning multiple times, we conclude that 
\begin{eqnarray*}
\sum\limits_{m_{k-1}\in\mathbb Z}\sum\limits_{m_{k-2}\in\mathbb Z}...\sum\limits_{m_{1}\in\mathbb Z}b_{n-m_{k-1}}b_{m_{k-1}-m_{k-2}}...b_{m_2-m_1}b_{m_1}\\&&\hspace{-5cm}=\sum\limits_{m_{k-1}\in\mathbb Z}b_{n-m_{k-1}}\left(\sum\limits_{m_{k-2}\in\mathbb Z}b_{m_{k-1}-m_{k-2}}\left(...\sum\limits_{m_{1}\in\mathbb Z}b_{m_2-m_1}b_{m_1}\right)\right),
\end{eqnarray*}
which proves the claim.
\end{remark}
\begin{proposition}[\cite{zrodlo5}]\label{rozdzielnosc}
Let $f,g,h\in\mathbb L$. Then
\begin{enumerate}
\item[(1)] $f\in\mathbb L (g)$, if and only if $g\in\mathbb L (f)$,
\item[(2)] $\mathbb L (g)\neq\emptyset$,
\item[(3)] $f\in\mathbb L (g)\Longrightarrow\alpha f\in\mathbb L (g)$ for every $\alpha\in\mathbb C$,
\item[(4)] $f,h\in\mathbb L (g)\Longrightarrow f+h\in\mathbb L (g)$.
\end{enumerate}
\end{proposition}
It is also easy to check that the product of formal Laurent series is commutative if it exists.

Since we will define and investigate properties of the composition of formal Laurent series in Section \ref{sec5}, we recall now for the convenience of the reader basic definitions and facts regarding composition of formal power series and formal Laurent series with formal power series. Recall that by $\mathbb X(\mathbb C)$ we denote the set of all formal power series over $\mathbb C$; because of the natural embeeding $\mathbb X(\mathbb C)\subset\mathbb L(\mathbb C)$, previous definitions of all operations on formal Laurent series are valid for formal power series (notice also that for every $g,f\in\mathbb X(\mathbb C)$, the product $gf\in\mathbb X(\mathbb C)$ is well-defined).
\begin{definition}[\cite{zrodlo8}]
Let $g=\sum\limits_{k=0}^{\infty}b_kz^k\in\mathbb X(\mathbb C)$. We define a subset $\mathbb X_g\subset\mathbb X(\mathbb C)$ as
\[
\mathbb X_g=\left\{f=\sum\limits_{k=0}^{\infty}a_kz^k\in\mathbb X(\mathbb C):\,\,\sum\limits_{n=0}^{\infty}b_na_k^{(n)}\mbox{ converges for all }k\in\mathbb Z_+\cup\{0\}\right\}.
\]
The mapping 
\[
T_g: \mathbb X_g\ni f=\sum\limits_{k=0}^{\infty}a_kz^k\mapsto\sum\limits_{k=0}^{\infty}\left(\sum\limits_{n=0}^{\infty}b_na_k^{(n)}\right)z^k\in\mathbb X(\mathbb C)
\]
is called the composition of $g$ and $f$; we denote $g\circ f:=T_g(f)$.
\end{definition}
In what follows, for every $f=\sum\limits_{k=0}^{\infty}a_kz^k\in\mathbb X(\mathbb C)$, deg$(f):=0$ if $f=0$, deg$(f):=+\infty$ if $a_k\neq 0$ for infinitely many $k$ and deg$(f):=\max\{k\in\mathbb Z_+\cup\{0\}:a_k\neq 0\}$ otherwise. We simply say that deg$(f)>0$ if deg$(f)\in\mathbb Z_+\cup\{+\infty\}$.
\begin{theorem}[\cite{zrodlo8}]\label{gencomp}
Let $f=\sum\limits_{k=0}^{\infty}a_kz^k,g=\sum\limits_{k=0}^{\infty}b_kz^k\in\mathbb X(\mathbb C)$, deg$(f)>0$. Then the composition $g\circ f$ exists, if and only if $\sum\limits_{n=k}^{\infty}\binom{n}{k}b_na_0^{n-k}$ converges for all $k\in\mathbb Z_+\cup\{0\}$ (equivalently: $g^{(k)}$ converges at $z=a_0$ for all $k\in\mathbb Z_+\cup\{0\}$, where $g^{(k)}$ denotes the $k$th formal derivative of $g$ treated as a classical power series).
\end{theorem}
The next two theorems state some of the most important properties of the composition of formal power series -- that is the Right Distributive Law and the Chain Rule, respectively.
\begin{theorem}[\cite{zrodlo9}, \cite{zrodlonew}, Thm. 5.6.2 and \cite{rem}, Remark 18]\label{rdlpower}
Let $f,g,h\in\mathbb X(\mathbb C)$, deg$(f)>0$. If $g\circ f$, $h\circ f$ exist, then $(gh)\circ f$ exists and $(gh)\circ f=(g\circ f)(h\circ f)$.
\end{theorem}
\begin{theorem}[\cite{zrodlo9}, \cite{zrodlonew}, Thm. 5.5.3]\label{crpower}
Let $f,g\in\mathbb X(\mathbb C)$, deg$(f)>0$. Then $g'\circ f$ exists, if and only if $g\circ f$ exists. Moreover, $(g\circ f)'=(g'\circ f)f'$ if $g\circ f$ exists.
\end{theorem}
\begin{remark}
Similarly to \cite{rem}, Remark 18, the assumption deg$(f)> 0$ must be included in the above theorem. Indeed, the proof of Theorem 5.5.3 in \cite{zrodlonew} is based on Lemma 5.5.2, which proof uses Theorem 5.4.1 (Thm. \ref{gencomp} in this paper). However, one of the assumptions of that theorem is deg$(f)> 0$, which is, unfortunately, omitted in Theorem 5.5.3 in \cite{zrodlonew}.
\end{remark}
Before introducing the next definition, it is worth mentioning that it is well-known (see e.g. \cite{zrodlonew}) that for $f=\sum\limits_{k=0}^{\infty}a_kz^k\in\mathbb X(\mathbb C)$, $a_0\neq 0$, there exists exactly one $g\in\mathbb X(\mathbb C)$ such that $gf=1$; we denote $g=f^{-1}$. Also, for $n\in\mathbb Z_+$, $f^{-n}:=(f^{-1})^n$ and $a_k^{(-n)}$ denote the coefficients of $f^{-n}$.
\begin{definition}[\cite{zrodlo5}]
Let $g=\sum\limits_{k\in\mathbb Z}b_kz^k\in\mathbb L (\mathbb C)$. Define
\[
\mathbb X_g=\left\{f=\sum\limits_{k=0}^{\infty}a_kz^k\in\mathbb X(\mathbb C):\,a_0\neq 0\mbox{ and }\sum\limits_{n\in\mathbb Z}b_na_k^{(n)}\mbox{ converges for all }k\in\mathbb Z_+\cup\{0\}\right\}.
\]
The mapping \[
T_g: \mathbb X_g\ni f=\sum\limits_{k=0}^{\infty}a_kz^k\mapsto\sum\limits_{k=0}^{\infty}\left(\sum\limits_{n\in\mathbb Z}b_na_k^{(n)}\right)z^k\in\mathbb X(\mathbb C)
\]
is called the composition of $g$ and $f$; we denote $g\circ f:=T_g(f)$.
\end{definition}

\begin{theorem}[\cite{zrodlo5}]
Let $g=\sum\limits_{k\in\mathbb Z}b_kz^k\in\mathbb L (\mathbb C)$ and $f=\sum\limits_{k=0}^{\infty}a_kz^k\in\mathbb X (\mathbb C)$, $a_0\neq 0$, deg$(f)>0$. Then $g\circ f$ exists, if and only if $\sum\limits_{n=k}^{\infty}\binom{n}{k}b_na_0^{n-k}$ and $\sum\limits_{n=k}^{\infty}\binom{n}{k}b_{-n}a_0^{k-n}$ converge for all $k\in\mathbb Z_+\cup\{0\}$.
\end{theorem}

\section{Basic properties of inverses of formal Laurent series}\label{sec3}

\begin{definition}\label{inverse}
We say that a formal Laurent series $g\in\mathbb L$ is an inverse of $f\in\mathbb L$, if and only if $fg=1$ (equivalently: $gf=1$). We denote the set of all inverses of a given $f\in\mathbb L$ as $R(f)$.
\end{definition}
It is well known that in the set $\mathbb X(\mathbb C)$ of formal power series, the inverse of a given series is unique if it exists; also, for $f=\sum\limits_{n=0}^{\infty}a_nz^n\in\mathbb X(\mathbb C)$, $f$ possesses an inverse in $\mathbb X(\mathbb C)$, if and only if $a_0\neq 0$ (see for example \cite{zrodlonew}, Thm. 1.1.8).  This case, however, seems to be much different for formal Laurent series -- as we will show, for $f\in\mathbb L$, exactly one of the following is true:
\begin{enumerate}
\item[(1)] $f$ has no inverses,
\item[(2)] $f$ has exactly one inverse,
\item[(3)] $f$ has uncountably many inverses.
\end{enumerate}
Let us introduce the following examples.

\begin{example}
Let $f=\sum\limits_{n\in\mathbb Z}a_nz^n$, where $a_n=a_0q^n$ for all $n\in\mathbb Z$, $a_0,q\in\mathbb C\setminus\{0\}$. Assume that $f$ possesses an inverse $g=\sum\limits_{n\in\mathbb Z}b_nz^n\in\mathbb L$. Then $0=\sum\limits_{m\in\mathbb Z}b_ma_{1-m}=q\sum\limits_{m=\in\mathbb Z}b_ma_{-m}=q$,
which is a contradiction since $q\neq 0$. This proves $R(f)=\emptyset$.
\end{example}

\begin{example}\label{uniquity}
Let $f=\sum\limits_{n\in\mathbb Z}a_nz^n$, where $a_n=1$ for $n\geq 0$ and $a_n=0$ for $n<0$. Assume that $g=\sum\limits_{n\in\mathbb Z}b_nz^n\in\mathbb L$ is an inverse of $f$. Then $\sum\limits_{m=0}^{\infty}b_{-m}=1$ and $\sum\limits_{m=0}^{\infty}b_{n-m}=0$ for $n\neq0$, so
\begin{eqnarray*}\label{coeff1}
b_n=\sum\limits_{m=0}^{\infty}b_{n-m}-\sum\limits_{m=0}^{\infty}b_{n-1-m}=\left\{\begin{array}{ll}
0, & n\in\mathbb Z\setminus\{0,1\},\\
-1, & n=1,\\
1, & n=0,\\
\end{array} \right.
\end{eqnarray*}
so $g=1-z$ (which is indeed an inverse of $f$, as can be checked explicitly). Therefore $f$ possesses exactly one inverse.
\end{example}

\begin{example}\label{laurentxx}
Let $f=\sum\limits_{n\in\mathbb Z}a_nz^n$, where $a_0=1$, $a_1=-1$ and $a_n=0$ for  $n\notin\left\{0,1\right\}$.
 A formal Laurent series $g=\sum\limits_{n\in\mathbb Z}b_nz^n$ is an inverse of $f$, if and only if
\[
\sum\limits_{m\in\mathbb Z}a_mb_{n-m}=b_n-b_{n-1}=\left\{\begin{array}{ll}
1, & n=0\\
0, & n\neq 0.\\
\end{array} \right.
\]
Therefore for every $c\in\mathbb C$, the formal Laurent series $\sum\limits_{n\in\mathbb Z}b_nz^n$, where $b_n=c$ if $n\geq 0$ and $b_n=c-1$ if $n<0$ is an inverse of $f$.
\end{example}
Examples \ref{uniquity} and \ref{laurentxx} show that the fact that $g$ is the only inverse of $f$ does not imply that $f$ is the only inverse of $g$. It can be also easily proved that a formal Laurent series cannot have a finite, but greater than 1, or a countable number of inverses. Indeed, we have the following
\begin{proposition}
Let $f\in\mathbb L$, $g_1,g_2\in R(f)$, $g_1\neq g_2$. Then $R(f)$ is uncountable. 
\end{proposition}
\begin{proof}
The claim is a direct result of the following observation: for every $k_1,k_2\in\mathbb C$, $k_2\neq -k_1$, $f\,\frac{k_1g_1+k_2g_2}{k_1+k_2}=\frac{k_1}{k_1+k_2}fg_1+\frac{k_2}{k_1+k_2}fg_2=1$. 
\end{proof}
\begin{remark}
Another two peculiar features of multiplication in $\mathbb L$ are worth noting, namely:
\begin{enumerate}
\item existence of zero divisors -- for example, $\left(\sum\limits_{n\in\mathbb Z}2^{-|n|}z^n-3\right)\sum\limits_{n\in\mathbb Z}z^n=0$ (\cite{zrodlo5}, Example 2.14),
\item non-associativity -- for example, let $f\in\mathbb L\setminus\{0\}$, $h_1,h_2\in R(f)$, $g:=h_1-h_2\neq 0$ (by Example \ref{laurentxx}, such $f$ exists); then $(h_1f)g=g\neq 0$ and $h_1(fg)=h_1(1-1)=0$.
\end{enumerate}
These two properties and the lack of general uniqueness of inverses in $\mathbb L$ turn out to be strongly connected. To show this, we will prove the non-associativity of multiplication in $\mathbb L$ implicitly, without constructing a particular example: assume that for every $f,g,h\in\mathbb L$, if $gh\in\mathbb L(f)$, then $fg\in\mathbb L(h)$ and $f(gh)=(fg)h$. Then 
\begin{enumerate}
    \item for any $f\in\mathbb L$, the existence of such $g\in\mathbb L\setminus\{0\}$ that $fg=0$ would imply $R(f)=\emptyset$ -- indeed, if $h\in R(f)$, we would have $g=g(fh)=(gf)h=0$;
    \item for every $f\in\mathbb L$, $R(f)$ would be either an empty set or a singleton -- if $g_1,g_2\in R(f)$, $g_1\neq g_2$, then $f(g_1-g_2)=0$ -- but that would mean $R(f)=\emptyset$ by (1).
\end{enumerate}
As a result, since there exist formal Laurent series possessing more than one inverse, the stated associativity assumption is false.
\end{remark}
One can, however, prove the associativity of multiplication of formal Laurent series satisfying some additional, but quite weak assumptions. For example, we have the following
\begin{proposition}\label{assoc}
Let $f=\sum\limits_{n\in\mathbb Z}a_nz^n,g=\sum\limits_{n\in\mathbb Z}b_nz^n,h=\sum\limits_{n\in\mathbb Z}c_nz^n\in\mathbb L$, $f\neq 0$. Define $|f|:=\sum\limits_{n\in\mathbb Z}|a_n|z^n$ (and $|g|$, $|h|$ analogously). Assume $|f|\in\mathbb L(|g|)$ and $|f||g|\in\mathbb L(|h|)$. Then
$(fg)h=f(gh)$ (in particular, all products in this formula exist).
\end{proposition}
\begin{proof}
Since $|f||g|\in\mathbb L(|h|)$, obviously $fg$ exists and the series $\sum\limits_{k\in\mathbb Z}\sum\limits_{s\in\mathbb Z}|a_sb_{k-s}c_{n-k}|$ converges for all $n\in\mathbb Z$. By the Rearrangement Theorem (see Remark \ref{rearr}), it follows that the series
\[
l_n:=\sum\limits_{k\in\mathbb Z}\sum\limits_{s\in\mathbb Z}a_sb_{k-s}c_{n-k}\mbox{ and }
r_n:=\sum\limits_{s\in\mathbb Z}\sum\limits_{k\in\mathbb Z}a_sb_{k-s}c_{n-k}
\]
are convergent and equal for all $n\in\mathbb Z$. Now, see that
\[
r_n=\sum\limits_{m\in\mathbb Z}\sum\limits_{k\in\mathbb Z}a_{n-m}b_{k-n+m}c_{n-k}=\sum\limits_{m\in\mathbb Z}\sum\limits_{t\in\mathbb Z}a_{n-m}b_{t}c_{m-t}
\]
(where we substituted $m=n-s$, $t=k-n+m$). Therefore, since $fg\in\mathbb L$, $(fg)h$ exists and, since $f\neq 0$, by a reasoning similar to Remark \ref{dlapopr}, $gh,f(gh)\in\mathbb L$. Also, $(fg)h=\sum\limits_{n\in\mathbb Z}l_nz^n=\sum\limits_{n\in\mathbb Z}r_nz^n=f(gh)$, which completes the proof.
\end{proof}
The lack of general associativity of multiplication and of general uniqueness of inverses in $\mathbb L$ makes the issue of composition of formal Laurent series more complex (according to our best knowledge, one can find in the literature the definition of composition of formal Laurent series and formal power series, and of composition of semi-Laurent series - see e.g. \cite{zrodlo5,biinfinite}). That's because it requires considering expressions of the form $f^{-n}:=(f^{-1})^n:=f^{-1}\cdot... \cdot f^{-1}$ ($n$ factors), and it is not possible to unambiguously assign to a given series $f$ its inverse $ f^{-1}$. We will, however, study this topic further in Section \ref{sec5}. 

\section{The existence and computation of inverses of formal Laurent series}\label{sec4}

In this section, we will consider the problem of existence and computation of inverses of formal Laurent series more explicitly -- we will provide some methods of constructing the set $R(f)$ for $f\in\mathbb L$ satisfying some minimal conditions required to employ the theory of infinite linear systems (see e.g. \cite{zrodlo1,zrodlo13,determinant,zrodlo2}). In particular, we will use the concept of the strictly particular solution to an infinite system of equations -- one of the particular solutions to the system, satisfying an infinite version of the Cramer rule (to see the full definition, the reader is referred e.g. to \cite{zrodlo1,zrodlo2}). It is known that an infinite system is consistent if and only if its strictly particular solution exists (also see e.g. \cite{zrodlo1,zrodlo2}).

Let $\f\in\mathbb L$. A series $g=\sum\limits_{n\in\mathbb Z}d_nz^n$ is an inverse of $f$, if and only if the following infinite system of linear equations is satisfied:
\begin{eqnarray}\label{ukladpierwszy}
\sum\limits_{m\in\mathbb Z}a_{n-m}d_m=\left\{\begin{array}{ll}
1,&n=0,\\
0,&n\neq 0.\\
\end{array}\right.
\end{eqnarray}
The index $n$ occurring above will be called the index of the equation. The methods for solving infinite systems mentioned above, however, require the equations and unknowns to be indexed by positive integers only. Therefore, the system (\ref{ukladpierwszy}) will first have to be transformed so that it meets this requirement. In the following subsections, we will present two approaches to this problem, obtaining two methods of constructing $R(f)$ (calculating inverses of $f$) under some assumptions. 
\begin{remark}
Let us note that, in accord with the theory of infinite systems presented in e.g. \cite{zrodlo1,zrodlo13,determinant,zrodlo2}, we define the determinant of an infinite matrix $A=[a_{i,j}]_{i,j\in\mathbb N}$ by the formula $|A|=\lim\limits_{n\rightarrow\infty}|A_n|$ (where $|A_n|$ is the determinant of the matrix $A_n=[a_{i,j}]_{i,j\in\{1,...,n\}}$) if this limit exists; otherwise we say that the determinant of $A$ does not exist. We call $|A_n|$ the $n$th principal minor of $A$. Also, we define the product of two infinite matrices $A=[a_{i,j}]_{i,j\in\mathbb N}$, $B=[b_{i,j}]_{i,j\in\mathbb N}$ by the formula $AB:=\left[\sum\limits_{k=1}^{\infty}a_{i,k}b_{k,j}\right]_{i,j\in\mathbb N}$ (provided all sums over $k$ converge) and say that $B$ is a both-side inverse of $A$, if and only if $AB=BA=I:=[\delta_{ij}]_{i,j\in\mathbb N}$.
\end{remark}

\subsection{Method I}\label{subsec4.1}
Let $\f\in\mathbb L$. Consider the infinite system 
\begin{eqnarray}\label{uklad}
W[f]\,y=b,
\end{eqnarray}
where $b=[b_1\,\,b_2\,\,...]^T:=[1\,\,0\,\,0\,\,...]^T$, $y=[y_1\,\,y_2\,\,...]^T$ are the unknowns and
\begin{eqnarray}\label{macierzwf}
W[f]=\left[ \begin{array}{lllll}
w_{1,1} & w_{1,2} &\ldots &w_{1,n} & \ldots\\
w_{2,1} & w_{2,2} &\ldots &w_{2,n} & \ldots\\
\vdots & \vdots & \ddots & \vdots &\ddots\\
w_{n,1} & w_{n,2} &\ldots & w_{n,n} & \ldots\\
\vdots &\vdots&\ddots&\vdots&\ddots\\
\end{array} \right] := \left[ \begin{array}{llllll}
a_0 & a_{-1} & a_1 &a_{-2} & a_2&\ldots\\
a_1 & a_0 &a_2 &a_{-1} &a_3& \ldots\\
a_{-1} & a_{-2} & a_0 & a_{-3} &a_1&\ldots\\
a_2 & a_1 &a_3 & a_0 & a_4&\ldots\\
a_{-2}&a_{-3}&a_{-1}&a_{-4}&a_0&\ldots\\
\vdots &\vdots&\vdots&\vdots&\vdots&\ddots\\
\end{array} \right],
\end{eqnarray}
that is $w_{i,j}=a_{(-1)^i\left\lfloor \frac{i}{2}\right\rfloor-(-1)^j\left\lfloor \frac{j}{2}\right\rfloor}$. We have the following
\begin{proposition}\label{ssystems}
\begin{enumerate}
\item if $(d_n)_{n\in\mathbb Z}$ is a solution to (\ref{ukladpierwszy}), then $(y_j)_{j\in\mathbb N}$, $y_j=d_{(-1)^j\lfloor\frac{j}{2}\rfloor}$ is a solution to (\ref{uklad});
\item if $(y_j)_{j\in\mathbb N}$ is a solution to (\ref{uklad}) and the series $\sum\limits_{m\in\mathbb Z}a_{n-m}d_m$, where
\begin{eqnarray}\label{odydod}
d_n=\left\{\begin{array}{ll}
y_{2n}, & n> 0\\
y_{1-2n}, & n\leq 0,\\
\end{array} \right. 
\end{eqnarray}
(that is $y_j=d_{(-1)^j\lfloor\frac{j}{2}\rfloor}$) converges for all $n\in\mathbb Z$, then  $(d_n)_{n\in\mathbb Z}$ is a solution to (\ref{ukladpierwszy}).
\end{enumerate}
\end{proposition}
\begin{proof}
Let $(d_n)_{n\in\mathbb Z}$ satisfy (\ref{ukladpierwszy}). Define $\mu_{n,s}=a_{n-(-1)^s\left\lfloor \frac{s}{2}\right\rfloor}$ for $n\in\mathbb Z$, $s\in\mathbb N$. It is easy to prove (by e.g. considering the Cauchy property for partial sums of the series $\sum\limits_{m=0}^{\infty}a_{n-m}d_m$, $\sum\limits_{m=0}^{\infty}a_{n+m}d_{-m}$ and $\sum\limits_{s=1}^{\infty}\mu_{n,s}y_s$) that since $\sum\limits_{m\in\mathbb Z}a_{n-m}d_m$ converges, $\sum\limits_{s=1}^{\infty}\mu_{n,s}y_s$ converges and is equal to $\sum\limits_{m\in\mathbb Z}a_{n-m}d_m=\delta_{n,0}$ for every $n\in\mathbb Z$. Now, changing the way of numbering the equations from $\mathbb Z$ to $\mathbb N$ by substituting $j=2n$ if $n>0$ and $j=1-2n$ if $n\leq 0$, we obtain (\ref{uklad}). The proof of (2) is similar.
\end{proof}

Now, let $\f$ satisfy the following conditions:
\begin{enumerate}
\item[(1)] $|W[f]|$ (the determinant of $W[f]$) exists and is different from $0$,
\item[(2)] all principal minors $|W[f]_k|$, $k\in\mathbb N$, of the matrix $W[f]$ are different from 0.
\end{enumerate}
\vskip .3cm
By [\cite{zrodlo1}, Thm. 2.1.], the matrix $W[f]$ can be written in the form
\begin{eqnarray}\label{gauss}
W[f]=\underbrace{\left[\begin{array}{llllll}
d_{1,1} & 0 &0&...&0&...\\
d_{2,1} &d_{2,2}&0&...&0&...\\
\vdots &\vdots&\ddots&\ddots&\vdots&\ddots\\
d_{n,1} & d_{n,2} &...&d_{n,n}&0&...\\
\vdots &\vdots&\ddots&\vdots&\vdots&\ddots\\
\end{array} \right] .}_{D}
\underbrace{\left[\begin{array}{llllll}
c_{1,1} & c_{1,2} &...&c_{1,n}&...\\
0 &c_{2,2}&...&c_{2,n}&...\\
\vdots &\vdots&\ddots&\vdots&\ddots\\
0 & 0 &... &c_{n,n}&...\\
\vdots &\vdots&\ddots&\vdots&\ddots\\
\end{array} \right] ,}_{C}
\end{eqnarray}
where (denoting $\sum\limits_{j=1}^{0}(...):=0$) the coefficients $d_{i,j}$, $c_{i,j}$ satisfy the recursive relations
\begin{eqnarray*}\label{bc}
d_{i,k}=\frac{w_{i,k}-\sum\limits_{j=1}^{k-1}d_{i,j}c_{j,k}}{c_{k,k}}\,\,(i\geq k),\quad c_{i,k}=\frac{w_{i,k}-\sum\limits_{j=1}^{i-1}d_{i,j}c_{j,k}}{d_{i,i}}\,\,(i< k),
\end{eqnarray*}
and the diagonal coefficients $c_{i,i}$, $i\in\mathbb N$ are arbitrary - let us then choose $c_{i,i}=1$ for all $i\in\mathbb N$ (see that for all $i\in\mathbb N$, $c_{i,i},d_{i,i}\neq 0$, because $|W[f]|\neq 0$). Now, by [\cite{zrodlo1}, Theorem 2.2, Corollary 2.2, Note 1], a both side inverse matrix $D^{-1}=[d^{(-1)}_{i,j}]_{i,j\in\mathbb N}$ of $D$ exists and is also triangular. It is then easy to check that, denoting $\beta_j:=d_{j,1}^{(-1)}$, we have the following recursion: 
\begin{eqnarray}\label{beta}
\left\{\begin{array}{l}
\beta_1=\frac{1}{d_{1,1}},\\
\beta_j=-\frac{1}{d_{j,j}}\sum\limits_{i=1}^{j-1}d_{j,i}\beta_i.\\
\end{array}\right.
\end{eqnarray}
The system (\ref{uklad}) is then equivalent [\cite{zrodlo1}, Cor. 2.2.] to the system $Cy=D^{-1}b$, that is
\begin{eqnarray}\label{gaussakrocej}
\sum\limits_{p=0}^{\infty}c_{j,j+p}y_{j+p}=\beta_j,\,\,j\in\mathbb N.
\end{eqnarray}

As mentioned before, the system (\ref{gaussakrocej}) is consistent (possesses a solution), if and only if its \textit{strictly particular solution} exists, that is, denoting
\begin{eqnarray}\label{doog}
A_p(j)=\left\{\begin{array}{ll}
\sum\limits_{k=0}^{p-1}(-1)^{p-1-k}c_{j+k,j+p}A_k(j),&p>0,\\
1,&p=0,\\
\end{array}\right.
\end{eqnarray}
the series
\begin{eqnarray}\label{rozwiazanie1}
B(j)=\sum\limits_{p=0}^{\infty}(-1)^pA_p(j)\beta_{j+p}
\end{eqnarray}
is convergent for all $j\in\mathbb N$ and $y_j=B(j)$ satisfies (\ref{gaussakrocej}). 
 It is worth mentioning that if (\ref{gaussakrocej}) is consistent, then [\cite{zrodlo1}, Thm. 3.8] $B(j)$ can be expressed as 
\begin{eqnarray}\label{rozwiazanie2}
B(j)=\frac{|C^{(j)}|}{|C|}=|C^{(j)}|\,\,\,\,\,(j\in\mathbb N),
\end{eqnarray}
where $|C|=\prod\limits_{i=1}^{\infty}c_{i,i}=1$ is the determinant of the matrix $C$ and $C^{(j)}$ is a matrix constructed by replacing the $j$th column of $C$ with the sequence $(\beta_i)_{i\in\mathbb N}$. See that we index the equations in (\ref{gaussakrocej}) by positive integers, unlike \cite{zrodlo1}, where indexes begin with 0, hence $j$ instead of $j+1$ in the above formula.\\
We therefore have the following
\begin{theorem}\label{jeden}
Let $f\in\mathbb L$ satisfy conditions \textit{(1)}-\textit{(2)}, that is $|W[f]|\in\mathbb C\setminus\{ 0\}$ and $|W[f]_k|\neq 0$ for all $k\in\mathbb N$. Then $R(f)\neq\emptyset$, if and only if $y_j:=B(j)$ calculated from (\ref{rozwiazanie2}) (or, equivalently, (\ref{rozwiazanie1})) exist and satisfy (\ref{gaussakrocej}), and the corresponding series $\sum\limits_{m\in\mathbb Z}d_ma_{n-m}$, $n\in\mathbb Z$, where $d_n$ are calculated from (\ref{odydod}), are convergent. Then $\sum\limits_{n\in\mathbb Z}d_nz^n$ is an inverse of $f$.
\end{theorem}
By Theorem \ref{jeden}, we obtain one inverse of $f\in\mathbb L$ satisfying \textit{(1)}-\textit{(2)} (if $R(f)\neq\emptyset$ of course). However, we know a formal Laurent series may have more than one inverse, so we now present a method of finding all inverses of $f$ satisfying \textit{(1)}-\textit{(2)}. We have the following simple
\begin{proposition}\label{all}
Let $f\in\mathbb L$ satisfy \textit{(1)}-\textit{(2)}, $R(f)\neq\emptyset$ and let $C$ denote the matrix from (\ref{gauss}). Then a formal Laurent series $g=\sum\limits_{n\in\mathbb Z}g_nz^n\in R(f)$, if and only if $g_n$ is of the form $g_n=d_n+d_n'$, where $d_n$ are calculated from Theorem \ref{jeden} and $y'_j:=d'_{(-1)^j\lfloor \frac{j}{2}\rfloor}$ satisfy a homogenous infinite system $C\,y'=0$.
\end{proposition}
\begin{proof}
Let $y^{(s)}=[y_1^{(s)}\,\,y_2^{(s)}\,...]^T$ be the strictly particular solution of (\ref{gaussakrocej}). See that $y=[y_1\,\,y_2\,...]^T\neq y^{(s)}$ is also a solution to (\ref{gaussakrocej}), if and only if, denoting $y'=y-y^{(s)}$,
$\sum\limits_{p=0}^{\infty}c_{j,j+p}y'_{j+p}=0$ ($j\in\mathbb N$). 
\end{proof}
Therefore obtaining the whole set $R(f)$ corresponds to solving the homogenous infinite system $C\,y'=0$. Although it is well-known that homogenous infinite systems are more complex to deal with that inhomogenous ones, there are necessary and sufficient conditions for their consistency and corresponding numerical algorithms of finding their solutions, see e.g. \cite{zrodlo1}, Thm 5.2 and \cite{fedorov2022}. In particular \cite{fedorov2022}, if the limits $x_j:=\lim\limits_{n\rightarrow\infty}\frac{(-1)^jA_{n-j}(j)}{A_n(0)}$ (where $A_p(j)$ were defined in (\ref{doog})) exist for all $j\in\mathbb N$, then $y_j'=cx_j$ satisfy $C\,y'=0$ for all $c\in\mathbb C$.
\begin{example}\label{exinf}
\begin{enumerate}
    \item Let $f=1-z$; since it is easy to check that $f$ satisfies conditions (1)-(2), we will now calculate $R(f)$ with the method presented above and compare the result with the one obtained in Example \ref{laurentxx}. In this case, formulas (\ref{gauss}), (\ref{beta}), (\ref{doog}), (\ref{rozwiazanie1}) give
    \[
    D=\begin{bmatrix}
1 & 0 & 0 & 0 & 0 & 0 & \cdots \\
-1 & 1 & 0 & 0 & 0 & 0 & \cdots \\
0 & 0 & 1 & 0 & 0 & 0 & \cdots \\
0 & -1 & -1 & 1 & 0 & 0 & \cdots \\
0 & 0 & 0 & 0 & 1 & 0 & \cdots \\
0 & 0 & 0 & -1 & -1 & 1 & \cdots \\
\vdots & \vdots & \vdots & \vdots & \vdots & \vdots & \ddots
\end{bmatrix},\quad C=\begin{bmatrix}
1 & 0 & -1 & 0 & 0 & 0 & \cdots \\
0 & 1 & -1 & 0 & 0 & 0 & \cdots \\
0 & 0 & 1 & 0 & -1 & 0 & \cdots \\
0 & 0 & 0 & 1 & -1 & 0 & \cdots \\
0 & 0 & 0 & 0 & 1 & 0 & \cdots \\
0 & 0 & 0 & 0 & 0 & 1 & \cdots \\
\vdots & \vdots & \vdots & \vdots & \vdots & \vdots & \ddots
\end{bmatrix},
    \]
$\beta_j=1$ for $j=1,2,4,6,8,...$, $\beta_j=0$ for $j=3,5,7,9,...$, $A_0(j)=1$ and $A_p(j)=0$ ($p>0$) for even $j$, $A_p(j)=1$ for even $p$ and odd $j$, $A_p(j)=0$ for odd $p$ and $j$ and, finally, $B(j)=\beta_j$ for every $j\in\mathbb N$. It is then easy to check that all conditions of Theorem \ref{jeden} are satisfied, and the obtained strictly particular solution corresponds to $1+z+z^2+...\in R(f)$. Now, the system $C\,y'=0$ in this example consists of equations $y'_1-y'_3=0$, $y'_2-y'_3=0$, $y'_3-y'_5=0$, $y'_4-y'_5=0$, ..., so its general solution is obviously $y'=[c\,\,c\,\,c\,\ \ldots]^T$, $c\in\mathbb C$. The obtained result is therefore in accord with Example \ref{laurentxx}.
    \item Let $f=\sum\limits_{n\in\mathbb Z}a_nz^n$, where
\[
a_n=\left\{\begin{array}{ll}
\sin n,&n>0,\\
\cos n,&n\leq 0.\\
\end{array}\right.
\]
It can be checked numerically that $f$ satisfies (1)-(2) since $|W[f]_k|\neq 0$ for $k\leq 50$ by an explicit calculation and a quick convergence of the sequence $(|W[f]_k|)$ to approximately $0.5$ was observed (up to accuracy of order $\sim 10^{-16}$ for $k\sim 50$) -- see upper-left plot of figure 1. Formulas (\ref{gauss}), (\ref{beta}), (\ref{rozwiazanie2}) were then used, with truncation to $k$ first rows and columns for increasing values of $k$ until satisfactory accuracy was reached. The results turned out to converge properly with increasing $k$, so we concluded that $R(f)\neq \emptyset$. More precisely, an approach similar to \cite{fedorov2022}, although with a more rigorous estimation of errors was used: for each $k$, the errors
\begin{eqnarray*}
\quad \quad\quad E_{1,k}=\max\limits_{1\leq j\leq k}|B^{(k+1)}(j)-B^{(k)}(j)|&,&\\E_{2,k}&=&\max\limits_{-\lfloor k/2\rfloor\leq n\leq \lfloor(k+1)/2\rfloor}\left|\sum\limits_{-\lfloor k/2\rfloor\leq m\leq \lfloor(k+1)/2\rfloor}d^{(k)}_ma_{n-m}-\delta_{n,0} \right|,
\end{eqnarray*}
were calculated (see lower part of figure 1), where $B^{(k)}(j)$ are numerical values of $B(j)$ calculated for truncation to $k$ columns and rows and $d^{(k)}_m$ are corresponding estimations of $d_m$ via (\ref{odydod}). Both $E_{1,k}$ and $E_{2,k}$ tend to $0$ for increasing $k$ (this is what we mean here by "proper convergence"); the values $d_m^{(k)}$ for the smallest value of $k$ such that $E_{1,k},E_{2,k}<10^{-10}$ were therefore taken as the numerical estimation of coefficients of $\sum\limits_{n\in\mathbb Z}d_nz^n\in R(f)$. They are also presented in figure 1. To determine whether $f$ possesses more than one inverse in this example (i.e. solving the homogenous system $C\,y'=0$), well-known (see e.g. \cite{fedorov2022}, sections 4-5), although more complex numerical procedures would be required; we omit this part here since the presented examples are mainly illustrative and detailed numerical analysis is not the subject of this paper. 

    \begin{figure}
    \includegraphics[width=0.99\textwidth]{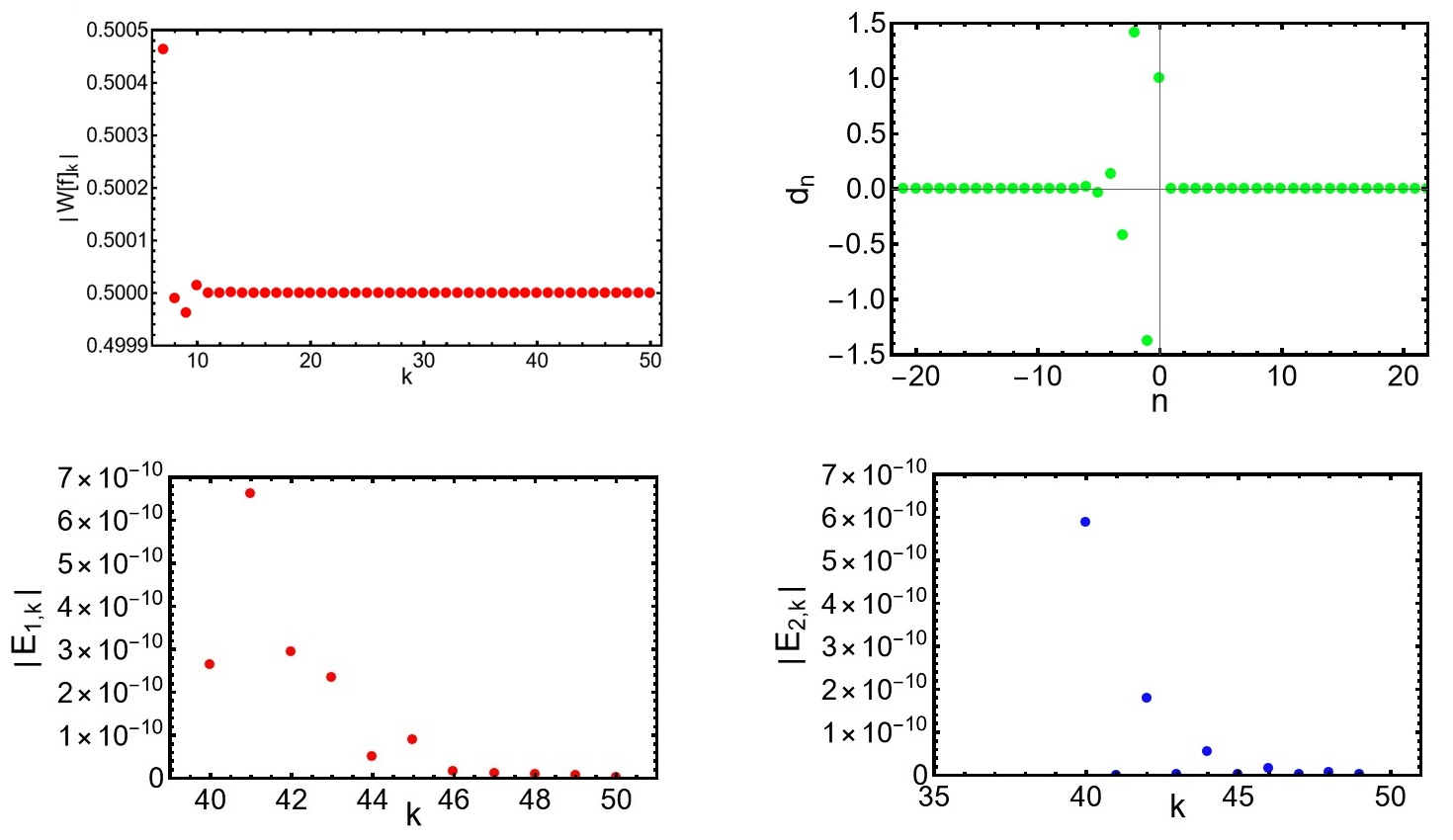}
    \caption{Upper-left plot: values of principal minors $|W[f]|_k$ of the matrix $W[f]$ (for $f$ from point (2) of Example \ref{exinf}) for various $k$; clearly a convergence towards approximately $0.5$ is visible. Lower plots: estimation of errors $E_{1,k}$, $E_{2,k}$ for various $k$; see that $k=44$ is the smallest $k$ for which $E_{1,k},E_{2,k}<10^{-10}$. Upper-right plot: values of $d_m^{(k)}$ for $k=44$.}
\end{figure}

    \item Let $f=\sum\limits_{n\in\mathbb Z}(|n|+1)z^n$. It can be easily checked numerically that the sequence $(|W[f]_k)_{k\in\mathbb N}$ diverges quickly to $-\infty$ (e.g. $|W[f]_{50}|\approx -1.436\cdot10^{16}$, $|W[f]_{100}|\approx-3.201\cdot 10^{31}$). However, if the method analogous to example (2) for finding the strictly particular solution is used, one obtains, for e.g. $N=50$, $d_{-24},d_{-23},...,d_{-2},d_2,...,d_{25}\approx 0$, $d_{-1},d_1\approx0.5$, $d_0\approx-1$ (with approximation errors not exceeding $2\cdot 10^{-15}$). Although the reliability of the numerical procedure could be questioned in this case, a simple analytical calculation can prove that indeed $\frac{1}{2z}-1+\frac12z\in R(f)$. This shows that even if the conditions (1)-(2) are not necessarily fulfilled, the method presented in this section may provide a useful intuition for analytical calculations (and/or a numerically convergent solution).
    
\end{enumerate}
\end{example}
\subsection{Method II}\label{subsec4.2}

We will now introduce a different approach to the problem (\ref{ukladpierwszy}), that is another method of finding the inverses of formal Laurent series satisfying some conditions (other than the ones considered in Section 4.1). As we will show, although it is less adapted to explicit numerical applications, it may be useful in some more general ("theoretical") considerations. Let us begin with the following
\begin{proposition}\label{sec2inf}
Let $f=\sum\limits_{n\in\mathbb Z}a_nz^n\in\mathbb L$. Then $g=\sum\limits_{n\in\mathbb Z}d_nz^n\in R(f)$, if and only if there exists a complex sequence $(s_j)_{j\in\mathbb N}$ such that sequences $x_i=d_{i-1}$ and $y_i=d_{-i}$ ($i\in\mathbb N$) satisfy the infinite systems
\begin{eqnarray}\label{uklady}
\sum\limits_{i=1}^{\infty}v_{i,j}x_i=s_j\,\,(j\in\mathbb N)\quad\mbox{ and }\quad \sum\limits_{i=1}^{\infty}u_{i,j}y_i=\left\{\begin{array}{ll}
1-s_j, &j=1,\\
-s_j,&j\in\mathbb N \setminus\{1\},\\
\end{array}\right.
\end{eqnarray}
where $v_{i,j}=a_{1-j+(-1)^i\left\lfloor \frac{i}{2}\right\rfloor}$, $u_{i,j}=a_{j+(-1)^i\left\lfloor \frac{i}{2}\right\rfloor}$ for $i,j\in\mathbb N$ -- that is, denoting $s=[s_1\,\,s_2\,\,...]^T$, $x=[x_1\,\,x_2\,\,...]^T$, $y=[y_1\,\,y_2\,\,...]^T$,
\begin{eqnarray*}
A_1[f]x=s,\qquad A_2[f]y=[1\,\,0\,\,0\,\,...]^T-s,
\end{eqnarray*}
where
\begin{eqnarray*}
A_1[f]=\left[\begin{array}{llllll}
a_0 &a_{-1}& a_{-2} &\ldots &a_{-n}&\ldots\\
a_1 &a_0 &a_{-1} &\ldots& a_{-n+1}&\ldots\\
a_{-1} &a_{-2} &a_{-3} &\ldots &a_{-n-1}&\ldots\\
a_2 &a_1 &a_0 &\ldots &a_{-n+2}&\ldots\\
a_{-2} &a_{-3} &a_{-4}& \ldots& a_{-n-2}&\ldots\\
\vdots &\vdots &\vdots &\ddots &\vdots &\ddots\\
\end{array} \right] ,\,\,
A_2[f]=\left[\begin{array}{llllll}
a_1 &a_2& a_3 &\ldots &a_n&\ldots\\
a_2 &a_3 &a_4 &\ldots& a_{n+1}&\ldots\\
a_0 &a_1 &a_2&\ldots &a_{n-1}&\ldots\\
a_3 &a_4 &a_5 &\ldots &a_{n+2}&\ldots\\
a_{-1} &a_0 &a_1& \ldots& a_{n-2}&\ldots\\
\vdots &\vdots &\vdots &\ddots &\vdots &\ddots\\
\end{array} \right] .
\end{eqnarray*}
\end{proposition}
\begin{proof}
    Assume $g\in R(f)$, that is $\sum\limits_{m\in\mathbb Z}d_ma_{n-m}=\delta_{n,0}$. In particular, $\sum\limits_{m=0}^{\infty}d_ma_{n-m}$ and $\sum\limits_{m=1}^{\infty}d_{-m}a_{n+m}$ are convergent, so there exists such a complex sequence $(u_n)_{n\in\mathbb Z}$ that \begin{eqnarray*}
\sum\limits_{m=0}^{\infty}a_{n-m}d_{m}=u_n\,\,(n\in\mathbb Z)\quad\mbox{and}\quad
\sum\limits_{m=1}^{\infty}a_{n+m}d_{-m}= \left\{\begin{array}{ll}
1-u_n, &n=0\\
-u_n,&n\neq 0\\
\end{array}\right. \,\,(n\in\mathbb Z),
\end{eqnarray*}
Rewriting the above system so that equations are indexed with positive integers only, we obtain
\[
\sum\limits_{m=0}^{\infty}a_{(-1)^j\left\lfloor j/2\right\rfloor-m}d_{m}=s_{j}\,\,(j\in\mathbb N)\quad\mbox{ and }\quad
\sum\limits_{m=1}^{\infty}a_{(-1)^j\left\lfloor j/2\right\rfloor+m}d_{-m}=\left\{\begin{array}{ll}
1-s_j, &j=1\\
-s_j,&j\neq 1\\
\end{array}\right.\,\,(j\in\mathbb N),
\]
where $s_j:=u_{(-1)^j\lfloor\frac{j}{2}\rfloor}$. The claim now follows directly from the definition of $x_i$, $y_i$, $v_{i,j}$ and $u_{i,j}$. The proof in the other direction is analogous.
\end{proof}
\begin{remark}
Denote the columns of $A_1[f]$ as $l_1$, $l_2$, ... and the columns of $A_2[f]$ as $m_1$, $m_2$, etc. Notice that $W[f]=[l_1\,\,l_2\,\,m_1\,\,l_3\,\,m_2\,\,l_4\,...]$.
\end{remark}
The above proposition shows that finding $R(f)$ for a given $f\in\mathbb L$ is equivalent to solving systems (\ref{uklady}) for all possible sequences $(s_j)_{j\in\mathbb N}$. In particular, $R(f)=\emptyset$, if and only if  for every $(s_j)$, at least one of the systems (\ref{uklady}) is inconsistent. Also, $f$ possesses a unique inverse, if and only if there exists exactly one sequence $(s_j)$ such that both systems (\ref{uklady}) are consistent, and both homogenous systems $A_1[f]x=0$, $A_2[f]y=0$ have no nontrivial solutions. This means we have the following
\begin{proposition}
A necessary condition for $f\in\mathbb L$ to possess exactly one inverse is that the kernels of the matrices $A_1[f]$ and $A_2[f]$ are both trivial.
\end{proposition}
\noindent Now, there are various possible ways to proceed from this point, for example:
\begin{enumerate}
\item solve both systems (\ref{uklady}) -- in particular, find all sequences $(s_j)$, for which both of them are consistent; if $A_{1,2}[f]$ have nonzero infinite determinants and nonzero principal minors, method analogous to that from Section 4.1 (LU decomposition and finding the strictly particular solution) can be used, 
\item find im$\,A_1[f]$ and determine the consistency of $A_2[f]y=[1\,0\,...]^T-s$ for all $s\in\mbox{im}\,A_1[f]$,
\item find im$\,A_1$, im$\,A_2$ and search for vectors $v_{1,2}$, $v_i\in$im$\,A_i[f]$ such that $v_{1}+v_2=[1\,0\,...]^T$.
\end{enumerate}
Although for a given $f\in\mathbb L$, these methods may be computationally difficult unless matrices $A_{1,2}[f]$ have a regular, simple form, they can be useful in more general investigations of the structure of $R(f)$ for some particular classes of formal Laurent series -- examples of which we are now going to provide.
\begin{example}
Let $f=a_0+a_1z+a_2z^2+...\in\mathbb L$, $a_0\neq 0$. We will find $R(f)$ (see that $f$ possesses a unique inverse in the set of formal power series but it may have more inverses in $\mathbb L$). It is easy to check that in this case, the systems (\ref{uklady}) reduce to $s_3=s_5=s_7=...=0$, $B_1x=[s_1\,\,s_2\,\,s_4\,\,s_6\,\,...]^T$, $B_2y=0$ and $B_3y=[1-s_1\,\,-s_2\,\,-s_4\,\,-s_6\,\,...]^T$, where
\begin{eqnarray*}
B_1=\left[\begin{array}{llll}
a_0 &0& 0 &\ldots\\
a_1 &a_0 &0 &\ldots\\
a_2 &a_1 &a_0 &\ldots \\
\vdots &\vdots &\vdots &\ddots \\
\end{array} \right] ,\,\,
B_2=\left[\begin{array}{llllll}
a_0 &a_1& a_2 &\ldots\\
0 &a_0 &a_1 &\ldots\\
0 &0 &a_0 &\ldots \\
\vdots &\vdots &\vdots &\ddots \\
\end{array} \right] ,\,\,
B_3=\left[\begin{array}{llll}
a_1 &a_2& a_3 &\ldots\\
a_2 &a_3 &a_4 &\ldots\\
a_3 &a_4 &a_5 &\ldots \\
\vdots &\vdots &\vdots &\ddots \\
\end{array} \right].
\end{eqnarray*}

Therefore $R(f)=\left\{\sum\limits_{n\in\mathbb Z}d_nz^n\in\mathbb L: B_2y=0,\,B_3y\mbox{ exists and }x=B_1^{-1}([1\,\,0\,\,...]^T-B_3y)\right\}$, where $x=[x_1\,\,x_2\,\,...]^T=[d_0\,\,d_1\,\,d_2\,\,...]^T$, $y=[y_1\,\,y_2\,\,...]^T=[d_{-1}\,\,d_{-2}\,\,...]^T$ (see that the trivial solution $y\equiv 0$ corresponds to the formal power series inverse of $f$). Let us emphasize that since $B_1$ is lower-triangular, its lower-triangular inverse $B^{-1}_1$ is easy to obtain by a recursive calculation of its coefficients:
\begin{eqnarray*}
\left\{\begin{array}{l}
B_{1\,m,m}^{(-1)}=\frac{1}{a_0}\\
B_{1\,n,m}^{(-1)}=-\frac{1}{a_0}\sum\limits_{i=m}^{n-1}a_{n-i}B_{1\,i,m}^{(-1)},\quad n>m.\\
\end{array}\right.
\end{eqnarray*}
Therefore constructing $R(f)$ reduces to solving the homogenous system $\sum\limits_{p=0}^{\infty}a_py_{j+p}=0$ ($j=1,2,...$), which is precisely the definition of a periodic linear infinite system, see e.g. \cite{periodic}. It has a nontrivial solution, if and only if its so called characteristic function $\phi(\mathcal X)=\sum\limits_{p=0}^{\infty}(-1)^pa_p\mathcal X^p$ has a nonzero radius of convergence and is equal to $0$ for some $\mathcal X$. If this condition is satisfied, the system's general solution is 
\[
y_j=(-1)^j\sum\limits_{k=1}^{\infty}S_k^j\sum\limits_{m=0}^{\nu_k-1}c_{m,k}j^m,
\]
where $S_k$ are zeros of $\phi(\mathcal X)$ of orders $\nu_k$ (of course the sum over $k$ is finite if $f$ has a finite number of zeros) and $c_{m,k}$ are constants such that the above series is convergent for each $j$ \cite{periodic}.\\
For instance, see that for $f=1+z+z^2+...$, $\phi(\mathcal X)=\sum\limits_{p=0}^{\infty}(-1)^p\mathcal X^p=\frac{1}{1+\mathcal X}$ ($|\mathcal X|<1$). Since $\phi$ has no zeros, $f$ has got only the formal power series inverse $1-z$, which is in agreement with Example \ref{uniquity}. On the other hand, for $g=1-z$, $\phi(\mathcal X)=1+\mathcal X$ has one zero: $S_1=-1$. Therefore $y_j=(-1)^j(-1)^jc_{0,1}j^0=c_{0,1}:=c$, $B_3y=[-c\,\,0\,\,0\,\,...]^T$ and $x_j=1+c$ for all $j\in\mathbb N$, which agrees with Example \ref{laurentxx}. Also, let us notice a more general fact -- if $f$ is a formal polynomial (a formal power series with finitely many nonzero terms), then the space of solutions to $B_2y=0$ is always finite-dimensional (since $\phi$ has finitely many zeros).
\end{example}
\begin{example}
Let $\f\in\mathbb L$ be a periodic formal Laurent series, that is there exists such $k\in\mathbb N$ that $a_{n+k}=a_n$ for all $n\in\mathbb Z$. Then for every $i,j\in\mathbb N$, $v_{i+2k,j}=v_{i,j}$ and  $u_{i+2k,j}=u_{i,j}$, so every $\alpha=[\alpha_1\,\,\alpha_2\,\,\alpha_3\,\,...]^T\in\mbox{im}A_1$ or im$\,A_2$ satisfies $\alpha_{i+2k}=\alpha_i$ for all $i\in\mathbb N$ (in particular dim$\,$im$\,A_{1,2}<\infty$). Therefore $[1\,\,0\,\,0\,\,...]^T$ cannot be written as $\alpha+\beta$, where $\alpha\in\mbox{im}\,A_1$, $\beta\in\mbox{im}\,A_2$. This proves $R(f)=\emptyset$ (see point (3) above).
\end{example}
\begin{example}
Let $\f\in\mathbb L$ be a symmetric formal Laurent series, that is satisfy $a_{-n}=a_n$ for all $n\in\mathbb Z$. Then an simple calculation shows that (\ref{uklady}) is equivalent to the pair of infinite systems
\[
A_2[f]x=[s_1\,\,s_3\,\,s_2\,\,s_5\,\,s_4\,...]^T-d_0[a_0\,\,a_1\,\,a_1\,\,a_2\,\,a_2\,...]^T,\qquad A_2[f]y=[1-s_1\,-s_2\,-s_3\,-s_4\,...]^T.
\]
This means finding $R(f)$ corresponds to finding all such sequences $(s_j)$ that $[s_1\,\,s_3\,\,s_2\,\,s_5\,\,s_4\,...]^T\in\mbox{im}\,A_2[f]\oplus \mbox{span}([a_0\,\,a_1\,\,a_1\,\,a_2\,\,a_2\,...]^T)$ and $[s_1\,\,s_2\,\,s_3\,\,s_4\,\,s_5\,...]^T\in\mbox{im}\,A_2[f]+\{[1\,\,0\,\,0\,\,...]^T\}$ (here $\oplus$ denotes the direct sum of linear spaces, while $+$ the Minkowski sum). In particular, if $g=\sum\limits_{n\in\mathbb Z}d_nz^n\in R(f)$ and $g^+f$ is symmetric like $f$, then (see proof of Prop. \ref{sec2inf}) $s_2=s_3$, $s_4=s_5$ etc., so $\left(\mbox{im}\,A_2[f]\oplus \mbox{span}([a_0\,\,a_1\,\,a_1\,\,a_2\,\,a_2\,...]^T)\right)\cap\left(\mbox{im}\,A_2[f]+\{[1\,\,0\,\,0\,\,...]^T\}\right)\neq\emptyset$ and therefore $[1\,\,0\,\,0\,\,...]^T-d_0[a_0\,\,a_1\,\,a_1\,\,a_2\,\,a_2\,...]^T\in\mbox{im}\,A_2[f]$.
\end{example}

\section{Composition of formal Laurent series}\label{sec5}
In this section, we are going to define the composition operation in the space of formal Laurent series, similarly to the well-known composition of formal power series (see Preliminaries). We will then focus on investigating its basic properties, e.g. whether it satisfies the Right Distributive Law (sec. 5.2) and the Chain Rule (sec. 5.3) similarly to formal power series.
\subsection{Definition and properties}\label{sec5.1}
Before defining properly the composition of two formal Laurent series, let us establish some necessary concepts and notations. Recall that for a formal Laurent series $f$, we define $f^0=1$, $f^1=f$ and $f^k=f^{k-1}f$ ($k>1$), provided $f^{k-1}$ exists and $f\in\mathbb L(f^{k-1})$. Let us now introduce the following
\begin{definition}\label{powers}
Let $f\in\mathbb L$, $R(f)\neq\emptyset$, $n\in\mathbb Z_+$. We define the $-n$th power of $f$ with respect to its inverse $f^{-1}\in R(f)$ as $f^{-n}:=(f^{-1})^n$, provided the right-hand side exists. Otherwise we say that $f^{-n}$ with respect to the inverse $f^{-1}$ does not exist; notice also that for different inverses of $f$, we may get a different value of $f^{-n}$. For $f=\sum\limits_{n\in\mathbb Z}a_nz^n$, $k\in\mathbb Z$, if $f^k$ exists, we denote $f^k=\sum\limits_{n\in\mathbb Z}a_n^{(k)}z^n$. 
\end{definition}
\begin{definition}
Let $\f\in\mathbb L$. We define deg$_+(f):=\sup\{n\in\mathbb Z_+\cup\{0\}:a_n\neq 0\}$ and deg$_-(f):=\sup\{n\in\mathbb Z_+:a_{-n}\neq 0\}$, where for the unification of notation we set $\sup\emptyset:=0$. Notice that deg$_+(f)$ (and similarly deg$_-(f)$) may be either a nonnegative integer or $+\infty$. We will simply write deg$_+(f)>0$ (deg$_-(f)>0$) if deg$_+(f)\in\mathbb Z_+\cup\{+\infty\}$ (deg$_-(f)\in\mathbb Z_+\cup\{+\infty\}$).
\end{definition}
Also, in what follows, we will denote $P_k(\mathbb L)=\{f\in\mathbb L:f^k\mbox{ exists}\}$ for every $k\in\mathbb Z_+\cup\{0\}$ and $P_{\infty}(\mathbb L)=\bigcap\limits_{k\in\mathbb N}P_k(\mathbb L)$. See that $\mathbb L=P_0(\mathbb L)=P_1(\mathbb L)\supset P_2(\mathbb L)\supset P_3(\mathbb L)\supset...$ . It is worth mentioning that although the multiplication of formal Laurent series involves infinite sums in general, the set $P_{\infty}(\mathbb L)$ is actually surprisingly "large", as it contains many important, broad classes of formal Laurent series -- e.g. formal power series, formal semi-Laurent and reversed semi-Laurent series (see e.g. \cite{zrodlonew}), the algebra $\mathbb L_1\equiv \ell^1(\mathbb Z)$ (see e.g. \cite{zrodlo5}) or Laurent expansions of analytic functions (again, see e.g. \cite{zrodlo5} for more details). 

Finally, we will use the notation $|f|:=\sum\limits_{n\in\mathbb Z}|a_n|z^n\in\mathbb L$ for $\f\in\mathbb L$. 

We are now going to define the composition of two formal Laurent series.
\begin{definition}\label{composition}
Let $\g\in\mathbb L$. Define
\begin{eqnarray*}
\,\mathbb X_g =\left\{\f\in P_{\mathrm{deg}_+(g)}(\mathbb L): \exists_{f^{-1}\in R(f)\cap P_{\mathrm{deg}_-(g)}(\mathbb L)}
\sum\limits_{k=-\mathrm{deg}_-(g)}^{\mathrm{deg}_+(g)}b_ka_n^{(k)}\mbox{ converges for all }n\in\mathbb Z\right\}
\end{eqnarray*}

if $g$ is not a formal power series (in the sense that deg$_-(g)>0$) and
\begin{eqnarray*}
\mathbb X_g =\left\{\f\in P_{\mathrm{deg}_+(g)}(\mathbb L):\sum\limits_{k=0}^{\mathrm{deg}_+(g)}b_ka_n^{(k)}\mbox{ converges for all }n\in\mathbb Z\right\}
\end{eqnarray*}
if $g$ is a formal power series (in the sense that deg$_-(g)=0$).\\
Now let $f\in\mathbb X_g$. If deg$_-(g)=0$, we define the composition of $g$ and $f$ as $g\circ f:=\sum\limits_{n\in\mathbb Z}\left(\sum\limits_{k=0}^{\mathrm{deg}_+(g)}b_ka_n^{(k)}\right)z^n$. If deg$_-(g)>0$, fix an inverse $f^{-1}\in R(f)$ such that the condition defining elements of $\mathbb X_g$ is satisfied; we define the composition $g\circ f$ with respect to the inverse $f^{-1}$ as $g\circ f:=\sum\limits_{n\in\mathbb Z}\left(\sum\limits_{k=-\mathrm{deg}_-(g)}^{\mathrm{deg}_+(g)}b_ka_n^{(k)}\right)z^n$.
\end{definition}
\begin{remark}\label{properties}
\begin{enumerate}
    \item If it does not cause any misunderstandings or if a given property holds for the composition $g\circ f$ with respect to any $f^{-1}\in R(f)$ satisfying the condition from the definition of $\mathbb X_g$, then we will simply call $g\circ f$ the composition of $g$ and $f$ instead of the composition of $g$ and $f$ with respect to $f^{-1}$.
    \item It is obvious that if $g\circ f,h\circ f$ exist (with respect to the same inverse of $f$ if deg$_-(g)$,deg$_-(h)>0$), then $f\in\mathbb X_{g+h}$ and $(g+h)\circ f=g\circ f+h\circ f$. Also, if $f\in\mathbb X_g$, then $f\in\mathbb X_{cg}$ and $(cg)\circ f=cg\circ f$ for all $c\in\mathbb C$.
    \item See that $g\circ f$ with respect to the inverse $f^{-1}$ of $f$ is equal to $\breve{g}\circ f^{-1}$ with respect to the inverse $f$ of $f^{-1}$ (of course if $g\circ f$ with respect to the inverse $f^{-1}$ exists).
    \item Assume $g'\circ f$ exists (with respect to a particular inverse $f^{-1}\in R(f)$ if deg$_-(g)<0$). Then, for all $n\in\mathbb Z$, $\sum\limits_{m\in\{-\mathrm{deg}_-(g),...,\mathrm{deg}_+(g)\}\setminus\{0\}}mb_{m}a_{n}^{(m-1)}$ converges, so (see e.g. \cite{zrodlo5}, Lemma 4.1) $\sum\limits_{m\in\{-\mathrm{deg}_-(g),...,\mathrm{deg}_+(g)\}\setminus\{0\}}b_{m}a_{n}^{(m-1)}$ converges as well, so $(z^{-1}g)\circ f$ exists.
\end{enumerate}
\end{remark}
As mentioned before, the lack of general uniqueness of inverses of formal Laurent series makes the issue of general composition $g\circ f$ ($g,f\in\mathbb L$) more complex. Fortunately, this ambiguity is not present if $g$ is a formal power series (deg$_-(g)=0$)-- in what follows we will focus mostly (but not only) on this case, which we believe is most useful for possible future applications. Let us mention, however, that for many important classes of formal Laurent series, one can overcome this non-uniqueness of $g\circ f$ even if deg$_-(g)>0$ -- for example, if $f$ is a so-called semi-formal Laurent series (that is deg$_-(f)>0$ is finite, see e.g. \cite{zrodlonew}), then $f\in P_{\infty}(\mathbb L)$ and $f$ always possesses exactly one inverse which is also semi-Laurent (the space of all semi-Laurent series is a field, see e.g. \cite{zrodlonew}) -- therefore compositions $g\circ f$ with respect to such inverses only can be considered.
\subsection{The Right Distributive Law}\label{subsec5.3}
We will now investigate whether a theorem analogous to Theorem \ref{rdlpower} is true for the composition of formal Laurent series. First, let us consider the following examples:
\begin{example}\label{rdl}
\begin{enumerate}
\item[(1)] Let $g=1-z$, $h=2z$, $f\in P_2(\mathbb L)$. It follows directly from Def. \ref{composition} that $g\circ f=1-f\in\mathbb L$, $h\circ f=2f\in\mathbb L$, $gh=2z-2z^2\in\mathbb L$, $(gh)\circ f=2f-2f^2\in\mathbb L$ and $(g\circ f)(h\circ f)=2f-2f^2\in\mathbb L$, so $(g\circ f)(h\circ f)=(gh)\circ f$.
\item[(2)] Let $g=\sum\limits_{n\in\mathbb Z}b_nz^n$, where $b_n=\frac{(-1)^{n}}{\sqrt{|n|}}$ for $n\neq 0$ and $b_0=0$. Let also $h=2g$ and $f=1$. It is obvious that the series $\sum\limits_{n\in\mathbb Z}b_n$ converges; denote $B:=\sum\limits_{n\in\mathbb Z}b_n$. By Def. \ref{composition}, $g\circ f=B\in\mathbb L$, $h\circ f=2B\in\mathbb L$, so $(g\circ f)(h\circ f)=2B^2\in\mathbb L$. However, see that
$\sum\limits_{m\in\mathbb Z}2b_mb_{-m}=\sum\limits_{m\in\mathbb Z\setminus\{0\}}\frac{2}{|m|}$, which of course diverges. Therefore $h\notin\mathbb L (g)$, so $(gh)\circ f$ does not exist.
\item[(3)] Let $g=z^l$, $h=z^m$ ($l,m\in\mathbb Z_+$) and $f\in P_{\max(l,m)}(\mathbb L)\setminus P_{l+m}(\mathbb L)$. Then $g\circ f,h\circ f\in\mathbb L$ and $gh\in\mathbb L$, but $(gh)\circ f$ does not exist. 
\item[(4)] Let $g=z^2$, $h=z^{-1}$ and $\f$, where $a_n=2$ for $n\geq 0$ and $a_n=1$ for $n<0$. See that $f\notin P_2(\mathbb L)$ and $1-z\in R(f)$. Therefore $h\circ f, (gh)\circ f\in\mathbb L$ (with respect to the inverse $1-z$ of $f$), but $g\circ f$ does not exist.
\item[(5)] Finally, let $g=h=z^2$, $f=\sum\limits_{n\in\mathbb Z}a_nz^n\in P_{2}(\mathbb L)$. We have $g\circ f=h\circ f=\sum\limits_{n\in\mathbb Z}a^{(2)}_nz^n\in\mathbb L$, where $a_n^{(2)}=\sum\limits_{m\in\mathbb Z}a_{n-m}a_{m}$. Now, $(gh)\circ f=\sum\limits_{n\in\mathbb Z}a^{(4)}_nz^n$ exists, if and only if the series
\begin{eqnarray}\label{pot1}
a^{(4)}_n=\sum\limits_{m_3\in\mathbb Z}\sum\limits_{m_2\in\mathbb Z}\sum\limits_{m_1\in\mathbb Z}a_{n-m_3}a_{m_3-m_2}a_{m_2-m_1}a_{m_1}
\end{eqnarray}
converges for all $n\in\mathbb Z$. Moreover, $(g\circ f)(h\circ f)$ exists, if and only if the series
\begin{eqnarray}\label{pot2}
\qquad\quad\sum\limits_{m_3\in\mathbb Z}\sum\limits_{m_2\in\mathbb Z}\sum\limits_{m_1\in\mathbb Z}a_{m_1}a_{m_3-m_1}a_{m_2}a_{n-m_3-m_2}=\sum\limits_{M_2\in\mathbb Z}\sum\limits_{M_3\in\mathbb Z}\sum\limits_{M_1\in\mathbb Z}a_{n-M_3}a_{M_3-M_2}a_{M_2-M_1}a_{M_1}
\end{eqnarray}
(where we denoted $M_1=m_1$, $M_2=m_3$ and substituted $M_3=m_2+m_3$) is convergent for all $n\in\mathbb Z$. However, see that the series (\ref{pot1}) and (\ref{pot2}) have the same terms but different order of summation, so one cannot state in general without a more careful analysis of the particular case that (\ref{pot1}) is convergent for all $n\in\mathbb Z$, if and only if (\ref{pot2}) is convergent for all $n\in\mathbb Z$ -- or that these series have equal sums even if they both converge.
\end{enumerate}
\end{example}
The above examples show that, unfortunately, the Right Distributive Law holds for some formal Laurent series but not in general (e.g. the existence of $(g\circ f)(h\circ f)$ does not necessarily imply the existence of $(gh)\circ f$). It can be proven, however, if some necessary additional assumptions are imposed -- to show that, let us begin with the following
\begin{lemma}\label{poww}
Let $\f\in\mathbb L$, $k\in\mathbb Z_+$ and assume $|f|^k$ exists. Then $f^0,...,f^k$ exist and $f^lf^{k-l}=f^k$ for all $l\in\{0,...,k\}$.
\end{lemma}
\begin{proof}
Since $|f|^k$ exists, the series $\sum\limits_{m_{k-1}\in\mathbb Z}\sum\limits_{m_{k-2}\in\mathbb Z}...\sum\limits_{m_{1}\in\mathbb Z}\left|a_{n-m_{k-1}}a_{m_{k-1}-m_{k-2}}...a_{m_2-m_1}a_{m_1}\right|$
is convergent for all $n\in\mathbb Z$ -- from that follows immediately that $f^k$ (and also $f^0,...,f^{k-1}$) exists since the series
\[
a_n^{(k)}=\sum\limits_{m_{k-1}\in\mathbb Z}\sum\limits_{m_{k-2}\in\mathbb Z}...\sum\limits_{m_{1}\in\mathbb Z}a_{n-m_{k-1}}a_{m_{k-1}-m_{k-2}}...a_{m_2-m_1}a_{m_1}
\]
is absolutely convergent (see Remark \ref{dlapopr}). Therefore (see Remark \ref{rearr}) the summation order can be rearranged and, for every $l\in\{2,...,k-2\}$ (for $l=0,1,k-1,k$ the claim is obvious), 
\[
a_n^{(k)}=\sum\limits_{m_{l}\in\mathbb Z}\sum\limits_{m_{k-1}\in\mathbb Z}...\sum\limits_{m_{l+1}\in\mathbb Z}\sum\limits_{m_{l-1}\in\mathbb Z}...\sum\limits_{m_{1}\in\mathbb Z}a_{n-m_{k-1}}a_{m_{k-1}-m_{k-2}}...a_{m_2-m_1}a_{m_1}.
\]
Denote $M_1=m_1$, ..., $M_{l-1}=m_{l-1}$, $M_0=m_l$, $M_l=m_{l+1}-m_l$, $M_{l+1}=m_{l+2}-m_l$, ..., $M_{k-2}=m_{k-1}-m_l$ -- then
\begin{eqnarray*}
a_n^{(k)}&=&\sum\limits_{M_0\in\mathbb Z}\sum\limits_{M_{k-2}\in\mathbb Z}...\sum\limits_{M_{1}\in\mathbb Z}a_{n-M_{k-2}-M_0}a_{M_{k-2}-M_{k-3}}...a_{M_{l+1}-M_l}a_{M_l}a_{M_0-M_{l-1}}a_{M_2-M_1}a_{M_1}\\
&=&\sum\limits_{M_0\in\mathbb Z}\left(\sum\limits_{M_{k-2}\in\mathbb Z}...\sum\limits_{M_{l}\in\mathbb Z}a_{n-M_0-M_{k-2}}a_{M_{k-2}-M_{k-3}}...a_{M_{l+1}-M_l}a_{M_l} \right)\cdot\\
&&\hspace{4.7cm}\cdot\left(\sum\limits_{M_{l-1}\in\mathbb Z}...\sum\limits_{M_{1}\in\mathbb Z} a_{M_0-M_{l-1}}a_{M_2-M_1}a_{M_1}\right)=\sum\limits_{M_0\in\mathbb Z}a_{n-M_0}^{(k-l)}a_{M_0}^{(l)}.\\
\end{eqnarray*}
\end{proof}

\begin{proposition}
Let $g=\sum\limits_{n=0}^{\infty}b_nz^n,h=\sum\limits_{n=0}^{\infty}c_nz^n\in\mathbb L$ (that is deg$_-(g)=\mbox{deg}_-(h)=0$) and let $f=\sum\limits_{n\in\mathbb Z}a_nz^n$ be a formal Laurent series such that the compositions $|g|\circ |f|$, $|h|\circ |f|$ exist and either: {\bf (a)} $(|g|\circ|f|)(|h|\circ |f|)$ exists or {\bf (b)} $(|g||h|)\circ |f|$ exists. Then $(g\circ f)(h\circ f),(gh)\circ f$ exist and $(g\circ f)(h\circ f)=(gh)\circ f$.
\end{proposition}
\begin{proof}
Let $k\in\left\{0,1,...,\max(\mbox{deg}_+(g),\mbox{deg}_+(h)\right\}$. Obviously $|f|^k$ exists; denote its coefficients as $|a_n|^{(k)}$. By formula (\ref{dopowers}) and Lemma \ref{poww}, $f^k$ exists as well and $\left|a_n^{(k)}\right|\leq|a_n|^{(k)}$ for every $n\in\mathbb Z$. Therefore the existence of $|g|\circ |f|$, $|h|\circ |f|$ implies that the series $\sum\limits_{k=0}^{\mathrm{deg}_+(g)}b_ka_n^{(k)}$, $\sum\limits_{k=0}^{\mathrm{deg}_+(h)}c_ka_n^{(k)}$ are absolutely convergent for all $n\in\mathbb Z$ (so $g\circ f$, $h\circ f$ exist) and, for all $n,m\in\mathbb Z$,
\begin{eqnarray}\label{pomm1}
\left(\sum\limits_{k=0}^{\mathrm{deg}_+(g)}b_ka_m^{(k)}\right)\left(\sum\limits_{k=0}^{\mathrm{deg}_+(h)}c_ka_{n-m}^{(k)}\right)=\sum\limits_{k=0}^{\infty}\sum\limits_{s=0}^kb_sc_{k-s}a_m^{(s)}a_{n-m}^{(k-s)},
\end{eqnarray}
where the right-hand side is also absolutely convergent. Let us clarify that we denote, e.g. if deg$_+(g)<+\infty$, $b_sa_m^{(s)}:=0$ for $s>\mbox{deg}_+(g)$ for simplicity (since $f^s$ may then not exist, the above expression could otherwise be ill-defined); we will use this notation in subsequent equations as well.\\
Let us first assume that condition {\bf (a)} holds. Then the following series is convergent for all $n\in\mathbb Z$: 
\[
\sum\limits_{m\in\mathbb Z}\left(\sum\limits_{k=0}^{\infty}|b_k||a_m|^{(k)}\right)\left(\sum\limits_{k=0}^{\infty}|c_k||a_{n-m}|^{(k)}\right)=\sum\limits_{m\in\mathbb Z}\sum\limits_{k=0}^{\infty}\sum\limits_{s=0}^k|b_s||c_{k-s}||a_m|^{(s)}|a_{n-m}|^{(k-s)}.
\]
See that by the Rearrangement Theorem (see Remark \ref{rearr}), 
$\sum\limits_{k=0}^{\infty}\sum\limits_{s=0}^k\sum\limits_{m\in\mathbb Z}|b_s||c_{k-s}||a_m|^{(s)}|a_{n-m}|^{(k-s)}$ converges for all $n\in\mathbb Z$; in particular, for $s_0:=\mbox{deg}_+(g)$, $k_0-s_0:=\mbox{deg}_+(h)$, $|b_{s_0}||c_{k_0-s_0}|\neq 0$, so $\sum\limits_{m\in\mathbb Z}|a_m|^{(s_0)}|a_{n-m}|^{(k_0-s_0)}$ converges -- and therefore by a reasoning analogous to the proof of Lemma \ref{dlapopr}, $|f|^{k_0}$ exists, so $|f|,f\in P(\mbox{deg}_+(g)+\mbox{deg}_+(h))$.\\
Now, since $|a_n^{(k)}|\leq |a_n|^{(k)}$, the series $r_n:=\sum\limits_{m\in\mathbb Z}\sum\limits_{k=0}^{\infty}\sum\limits_{s=0}^kb_sc_{k-s}a_m^{(s)}a_{n-m}^{(k-s)}$ is absolutely convergent for all $n\in\mathbb Z$. Therefore: (1) by eq. (\ref{pomm1}), $\sum\limits_{m\in\mathbb Z}\left(\sum\limits_{k=0}^{\mathrm{deg}_+(g)}b_ka_m^{(k)}\right)\left(\sum\limits_{k=0}^{\mathrm{deg}_+(h)}c_ka_{n-m}^{(k)}\right)$ converges for all $n\in\mathbb Z$, so $(g\circ f)(h\circ f)$ exists, (2) by the Rearrangement Theorem same as before and Lemma \ref{poww}, $r_n=\sum\limits_{k=0}^{\infty}\sum\limits_{s=0}^k\sum\limits_{m\in\mathbb Z}b_sc_{k-s}a_m^{(s)}a_{n-m}^{(k-s)}=\sum\limits_{k=0}^{\mathrm{deg}_+(g)+\mathrm{deg}_+(h)}\left(\sum\limits_{s=0}^kb_sc_{k-s}\right)a_n^{(k)}$, which proves the existence of $(gh)\circ f$ and the equality $(g\circ f)(h\circ f)=(gh)\circ f$.\\
Now assume {\bf (b)} holds. Then obviously $|f|,f\in P_{\mathrm{deg}_+(g)+\mathrm{deg}_+(h)}(\mathbb L)$. For every $n\in\mathbb Z$, the series (see Lemma \ref{dlapopr})
\[
\sum\limits_{k=0}^{\mathrm{deg}_+(g)+\mathrm{deg}_+(h)}\left(\sum\limits_{s=0}^k|b_s||c_{k-s}|\right)|a_n|^{(k)}=\sum\limits_{k=0}^{\mathrm{deg}_+(g)+\mathrm{deg}_+(h)}\sum\limits_{s=0}^k\sum\limits_{m\in\mathbb Z}|b_s||c_{k-s}||a_m|^{(s)}|a_{n-m}|^{(k-s)}
\]
converges. Analogously to {\bf (a)}, it can be proven that $(gh)\circ f$ exists, the series 
\[
\sum\limits_{m\in\mathbb Z}\sum\limits_{k=0}^{\infty}\sum\limits_{s=0}^kb_sc_{k-s}a_m^{(s)}a_{n-m}^{(k-s)}=\sum\limits_{m\in\mathbb Z}\left(\sum\limits_{k=0}^{\mathrm{deg}_+(g)}b_ka_m^{(k)}\right)\left(\sum\limits_{k=0}^{\mathrm{deg}_+(h)}c_ka_{n-m}^{(k)}\right)
\]
converges for all $n\in\mathbb Z$, so $(g\circ f)(h\circ f)$ exists and $(gh)\circ f=(g\circ f)(h\circ f)$, which completes the proof.
\end{proof}
\begin{corollary}
Let $f=\sum\limits_{n\in\mathbb Z}a_nz^n,g=\sum\limits_{n=0}^{\infty}b_nz^n,h=\sum\limits_{n=0}^{\infty}c_nz^n\in\mathbb L$. Assume $a_n,b_n,c_n\in\mathbb R_+\cup\{0\}\subset\mathbb C$ for all $n$ and $g\circ f$, $h\circ f$ exist. Then $(g\circ f)(h\circ f)$ exists, if and only $(gh)\circ f$ exists; moreover, $(g\circ f)(h\circ f)=(gh)\circ f$ if either side is well-defined.
\end{corollary}

\begin{remark}
The situation seems to be more complicated if we allow deg$_-(g),\mbox{deg}_-(h)>0$. Indeed, let for example $f\in P_2(\mathbb L)$ and let $f^{-1}:=\sum\limits_{n\in\mathbb Z}j_nz^n\in R(f)$. Then it is easy to check that for every $n\in\mathbb Z$, $a_n=\sum\limits_{m\in\mathbb Z}\sum\limits_{k\in\mathbb Z}a_ma_{k-m}j_{n-k}$, while $f^2f^{-1}$ exists, if and only if the following series is convergent for all $n\in\mathbb Z$: $\sum\limits_{k\in\mathbb Z}\sum\limits_{m\in\mathbb Z}a_ma_{k-m}j_{n-k}$. Therefore, without imposing any additional absolute convergence conditions, one cannot a priori assume that $f^2f^{-1}=f$ (even if e.g. $|f|\in P_2(\mathbb L)$ -- it is not even necessarily true that $|f^{-1}||f|$ exists).
\end{remark}

\begin{remark}
Let $f,g\in\mathbb L$ and let $g^{-1}$ be an inverse of $g$. Assume $g^{-1}\circ f$, $g\circ f$ exist. See that $g^{-1}\circ f$ is an inverse of $g\circ f$, if and only if $g,g^{-1},f$ satisfy the Right Distributive Law, that is $(g\circ f)(g^{-1}\circ f)=1=(gg^{-1})\circ f$.
\end{remark}

\subsection{The Chain Rule}\label{subsec5.2}
We will now move on to the problem whether a theorem analogous to Theorem \ref{crpower} is true for the composition of formal Laurent series. Similarly to the previous subsection, let us begin by considering some examples:
\begin{example}\label{cr}
\begin{enumerate}
\item[(1)] Let $g=1-z$ and $f\in\mathbb L$. We have $g'=-1$, so $g'\circ f=-1\in\mathbb L$ and $(g'\circ f)f'=-f'\in\mathbb L$. Moreover, $(g\circ f)'=(1-g)'=-g'\in\mathbb L$, so $(g\circ f)'=(g'\circ f)f'$.
\item[(2)] Let $g=\sum\limits_{n\in\mathbb Z\setminus\{0\}}\frac{1}{n^2}z^n$ and $f=1$. Since $\sum\limits_{m\in\mathbb Z\setminus\{0\}}\frac{1}{m^2}$ converges, $g\circ f$ exists (so $(g\circ f)'$ exists). However, $g'\circ f$ does not exist since $\sum\limits_{m\in\mathbb Z\setminus\{-1\}}\frac{1}{m+1}$ diverges, so $(g'\circ f)f'$ does not exist (even though $g'=0$).
\item[(3)] Let $g=z^2$ and $f=\sum\limits_{n\in\mathbb Z\setminus\{0\}}\frac{1}{n}z^n$. See that $f\in P_2(\mathbb L)$ (so $g\circ f$ exists) since $\sum\limits_{m\in\mathbb Z\setminus\{0,n\}}\frac{1}{m}\frac{1}{n-m}$ converges for all $n\in\mathbb Z$. Also, $g'\circ f=2f\in\mathbb L$. However, see that $ff'$ (and therefore $(g'\circ f)f'$) does not exist, since the series $\sum\limits_{m\in\mathbb Z\setminus\{0,n+1\}}\frac{1}{m}(n-m+1)\frac{1}{n-m+1}$ diverges for all $n\in\mathbb Z$.
\end{enumerate}
\end{example}
The above examples show that the Chain Rule, like the Right Distributive Law, holds for some but not all formal Laurent series; it does not even always hold for $g=z^k$, $k\geq 2$ - indeed, example (3) shows that one can have $f\in P_k(\mathbb L)$ but $f^{k-1}\notin \mathbb L(f')$. Again, however, it holds for formal Laurent series satisfying some additional assumptions. To prove that, let us first introduce the following
\begin{lemma}\label{docr}
Let $\f\in\mathbb L$, $\mathcal N\in\mathbb Z_+$, $N>1$; assume $|f|\in P_{\mathcal N-1}(\mathbb L)$ and $|f|,...,|f|^{\mathcal N-1}\in\mathbb L(|f|')$. Then for all $N\in\{1,...,\mathcal N\}$, $(f^N)'=Nf^{N-1}f'$.
\end{lemma}
\begin{proof}
Let $N\in\{1,...,\mathcal N-1\}$. Since $|f|^N$ exists, then, by Lemma \ref{poww} and formula (\ref{dopowers}), $f^N$ exists and $\left|a_n^{(N)}\right|\leq |a_n|^{(N)}$ for all $n\in\mathbb Z$. Now, the existence of $|f|^{N}|f|'$ means that $\sum\limits_{m=0}^{\infty}(m+1)|a_{m+1}||a_{n-m}|^{(N)}$, $\sum\limits_{m=1}^{\infty}(-m+1)|a_{-m+1}||a_{n+m}|^{(N)}$ converge for all $n\in\mathbb Z$; since the corresponding series with $|m+1|$, $|-m+1|$ instead of $m+1$, $-m+1$ obviously converge as well, it follows that $\sum\limits_{m\in\mathbb Z}|m+1||a_{m+1}||a_{n-m}|^{(N)}$ converges for all $n$ -- and therefore, by the inequality $\left|a_n^{(N)}\right|\leq |a_n|^{(N)}$, $f^N\in\mathbb L(f')$. In particular, $f^{\mathcal N-1}\in\mathbb L(f')$, so by \cite{zrodlo5}, Theorem 4.2 (iii), $f^{\mathcal N}$ exists. Also, arguing similarly as in the proof of Proposition \ref{assoc} and the beginning of this proof, it is easy to show that $f^{N-1}f'=(ff^{N-2})f'=f(f^{N-2}f')$ etc.\\
Therefore, using \cite{zrodlo5}, Theorem 4.2 (v), we have, for all $N\in\{1,...,\mathcal N\}$,
\begin{eqnarray*}
(f^N)'&=&f^{N-1}f'+f(f^{N-1})'=f^{N-1}f'+f(f^{N-2}f'+f(f^{N-2})')=\,...\\
&=&(...(((2ff')f+f^2f')f+f^3f')f+...+f^{N-2}f')f+f^{N-1}f'\\
&=&(...(f(3f^2f')+f^3f')f+...+f^{N-2}f')f+f^{N-1}f'\\
&=&(...(4f^3f')f+...+f^{N-2}f')f+f^{N-1}f'=...=Nf^{N-1}f',
\end{eqnarray*}
which completes the proof.
\end{proof}
\begin{proposition}
Let $g=\sum\limits_{n=0}^{\infty}b_nz^n,\f\in\mathbb L$ and assume $(|g|'\circ|f|)|f|'$ exists. Then $(g'\circ f)f',\,\,(g\circ f)'$ exist and are equal.
\end{proposition}
\begin{proof}
If deg$_+(g)=0$, the claim is obvious; let us therefore assume deg$_+(g)>0$.\\
Since $(|g|'\circ|f|)|f|'$ exists, $\sum\limits_{k\Z}(n-k+1)|a_{n-k+1}|\sum\limits_{m=1}^{\mathrm{deg}_+(g)}m|b_m||a_k|^{(m-1)}$ is convergent for all $n\in\mathbb Z$ -- therefore so are the series $\sum\limits_{k=-\infty}^n|n-k+1||a_{n-k+1}|\sum\limits_{m=0}^{\mathrm{deg}_+(g)}m|b_m||a_k|^{(m-1)}$ and also  $\sum\limits_{k=n+1}^{\infty}|n-k+1||a_{n-k+1}|\sum\limits_{m=0}^{\mathrm{deg}_+(g)}m|b_m||a_k|^{(m-1)}$ (where we could add the absolute value to $n-k+1$ due to its constant sign on each of the distinguished subsets of $\mathbb Z$). As a result, $\sum\limits_{k\Z}|n-k+1||a_{n-k+1}|\sum\limits_{m=1}^{\mathrm{deg}_+(g)}m|b_m||a_k|^{(m-1)}$ converges for $n\in\mathbb Z$, so by Lemma \ref{poww} and the inequality $|a_n^{(k)}|\leq |a_n|^{(k)}$, $r_n:=\sum\limits_{k\Z}\sum\limits_{m=1}^{\mathrm{deg}_+(g)}(n-k+1)a_{n-k+1}mb_ma_k^{(m-1)}$ (and obviously all the series $\sum\limits_{m=1}^{\mathrm{deg}_+(g)}mb_ma_k^{(m-1)}$) is absolutely convergent. This proves the existence of $(g'\circ f)f'$; also, by the Rearrangement Theorem (see Remark \ref{rearr}), $r_n=\sum\limits_{m=1}^{\mathrm{deg}_+(g)}\sum\limits_{k\Z}(n-k+1)a_{n-k+1}mb_ma_k^{(m-1)}$. Now, see that: (1) if $\mathrm{deg}_+(g)<+\infty$, then since $\mathrm{deg}_+(g)b_{\mathrm{deg}_+(g)}\neq 0$, obviously $\sum\limits_{k\Z}(n-k+1)|a_{n-k+1}||a_k|^{(\mathrm{deg}_+(g)-1)}$ converges, so $|f|^{\mathrm{deg}_+(g)-1}\in\mathbb L(|f|')$; by a reasoning analogous to the beginning of the proof and the proof of Prop. \ref{assoc}, if $\mathrm{deg}_+(g)\geq 2$, then $|f|^{\mathrm{deg}_+(g)-1}|f|'=|f|(|f|^{\mathrm{deg}_+(g)-2}|f|')$ -- in particular, $|f|^{\mathrm{deg}_+(g)-2}|f|'$ exists and, analogously, $|f|,...,|f|^{\mathrm{deg}_+(g)-1}\in\mathbb L(|f|')$ (2) if $\mathrm{deg}_+(g)=+\infty$, then by a similar reasoning we conclude that $|f|\in P_{\infty}(\mathbb L)$ and all natural powers of $|f|$ belong to $L(|f|')$. Therefore, in either case ((1), (2)), we conclude from Lemma \ref{docr} that for all $1\leq N\leq\mathrm{deg}_+(g)$ ($<\mathrm{deg}_+(g)$ if $\mathrm{deg}_+(g)=+\infty$), $(f^N)'=Nf^{N-1}f'$ (notice that the proof of Lemma \ref{docr} ensures that both sides are well-defined), which means $(n+1)a_{n+1}^{(N)}=N\sum\limits_{k\in\mathbb Z}a_k^{(N-1)}(n-k+1)a_{n-k+1}$. Therefore $r_n=(n+1)\sum\limits_{m=1}^{\mathrm{deg}_+(g)}b_ma_{n+1}^{(m)}=(n+1)\sum\limits_{m=0}^{\mathrm{deg}_+(g)}b_ma_{n+1}^{(m)}$ (since $(n+1)a_{n+1}^{(0)}=(n+1)\delta_{n+1,0}=0$). Therefore $g\circ f$ exists and $(g\circ f)'=(g'\circ f)f'$, which completes the proof.
\end{proof}
\vskip.5cm
Finally, let us emphasize that the problem of finding general necessary or/and sufficient conditions for the existence of composition of any two given formal Laurent series remains open.

\section*{Acknowledgements}

I would like to thank prof. Piotr Ma\'{c}kowiak for his comments, which have helped me improve the section "The existence and computation of inverses of formal Laurent series".

\end{document}